\documentclass[11pt]{article}
\usepackage{latexsym,amssymb}
\usepackage{amsmath}
\usepackage[mathscr]{eucal}
\usepackage{color}
\usepackage[abbrev]{amsrefs}

\newtheorem{Theorem}{\bf Theorem}[section]
\newtheorem{Lemma}{\bf Lemma}[section]
\newtheorem{Proposition}{\bf Proposition}[section]
\newtheorem{Corollary}{\bf Corollary}[section]
\newtheorem{Remark}{\bf Remark}[section]
\newtheorem{Example}{\bf Example}[section]
\newtheorem{Definition}{\bf Definition}[section]

\newenvironment{theorem}{\begin{Theorem}$\!\!\!$}{\end{Theorem}}
\newenvironment{lemma}{\begin{Lemma}$\!\!\!$}{\end{Lemma}}
\newenvironment{proposition}{\begin{Proposition}$\!\!\!$}{\end{Proposition}}

\newenvironment{remark}{\begin{Remark}$\!\!\!$}{\end{Remark}}

\newenvironment{definition}{\begin{Definition}$\!\!\!$}{\end{Definition}}

\numberwithin{equation}{section}

\numberwithin{equation}{section}

\newcommand{\dee}{{\rm{d}}}

\def\XXint#1#2#3{{\setbox0=\hbox{$#1{#2#3}{\int}$}
\vcenter{\hbox{$#2#3$}}\kern-.5\wd0}}

\allowdisplaybreaks
\begin{document}

\title{A one-dimensional Stefan problem for the heat equation\\ 
with a nonlinear boundary condition}
\author{
\qquad\\
Kensho Araya and Kazuhiro Ishige\\ \\
Graduate School of Mathematical Sciences, The University of Tokyo,\\
\qquad\,\,\, 3-8-1 Komaba, Meguro-ku, Tokyo 153-8914, Japan. }
\date{}
\maketitle
\begin{abstract}
We study the one-dimensional one-phase Stefan problem for the heat equation with a nonlinear boundary condition.
We show that all solutions fall into one of three distinct types: 
global-in-time solutions with exponential decay, global-in-time solutions with non-exponential decay, and finite-time blow-up solutions.
The classification depends on the size of the initial function. 
Furthermore, we describe the behavior of solutions at the blow-up time.
\end{abstract}
\vspace{35pt}
\smallskip
\noindent
E-mail: {\tt kensho@ms.u-tokyo.ac.jp} (K.A.),  {\tt ishige@ms.u-tokyo.ac.jp} (K.I.)\\

\vspace{3pt}
\newpage
%%%%%%%%%%%%%%%%%%%%%%%%%%%%%%
%%%%%%%%%%%%%%%%%%%%%%%%%%%%%%
\section{Introduction}
%%%%%%%%%%%%%%%%%%%%%%%%%%%%%%
%%%%%%%%%%%%%%%%%%%%%%%%%%%%%%
Consider the one-dimensional one-phase Stefan problem for the heat equation with a nonlinear boundary condition:
\[
\tag{P}
\label{SP}
\left\{
\begin{array}{ll}
\partial_t u=\partial_x^2 u, & (t,x)\in D_s(T),\vspace{3pt}\\
-\partial_xu(t,0)=u(t,0)^p,\qquad & t\in(0,T], \vspace{3pt}\\
u(t,s(t))=0, &  t\in (0,T],\vspace{3pt}\\
s'(t)=-\partial_xu(t,s(t)), &  t\in (0,T], \vspace{3pt}\\
u(0,x)=\varphi(x), & x\in[0,s_0],\vspace{3pt}\\
s(0)=s_0,
\end{array}
\right.
\]
where $\partial_t:=\partial/\partial t$, $\partial_x:=\partial/\partial x$, $T\in(0,\infty]$, $p>1$, and $(s_0,\varphi)\in W$. 
Here 
\begin{align}
\label{eq:1.1}
 & D_s(T):=\bigcup_{t\in(0,T]}\{t\}\times(0,s(t)),\\
 \nonumber
 & W:=\{(a,\phi)\,:\,a\in(0,\infty),\,\,\,\phi\in \mbox{Lip}([0,a])\setminus\{0\}\mbox{ with $\phi\ge 0$ on $[0,a]$ and $\phi(a)=0$}\}.
 \end{align}
We adopt the convention that $(0,T]=(0,T)$ if $T=\infty$.
This model describes heat diffusion in water under two types of boundary conditions:
a nonlinear reaction condition at $x=0$, representing an exothermic reaction,
and a Stefan-type condition at the moving boundary $x=s(t)$, representing the phase transition from ice to water.
Stefan problems with nonlinear boundary conditions have been investigated in many papers; 
see, for example,  \cites{DKS, Ke1, Ke2, Ke3, K, NP, P}. 

Given any $(s_0,\varphi)\in W$, problem~\eqref{SP} admits a unique solution $(s,u)$ on $(0,T]$ for some $T\in(0,\infty]$ 
(see Section~2). 
The maximal existence time of the solution is denoted by $T_m$, i.e., $T_m$ is the supremum of such $T$.
If $T_m<\infty$, then 
\begin{equation}
\label{eq:1.2}
\limsup_{t\nearrow T_m} u(t,0)=\infty
\end{equation}
(see Proposition~\ref{Proposition:2.4}), 
and $T_m$ is called the blow-up time of the solution $(s,u)$.
Furthermore, it follows from the maximum principle and Hopf's lemma that 
\begin{equation}
\label{eq:1.3}
u>0\quad\mbox{in}\quad D_s(T),\qquad -\partial_x u(t,s(t))=s'(t)>0\quad\mbox{for}\quad t\in (0,T], 
\end{equation}
for $T\in(0,T_m)$. 
In this paper, we study the behavior of solutions $(s,u)$ to problem~\eqref{SP} and classify them into three distinct types:
\begin{itemize}
  \item[(a)]  
  Global-in-time solutions with exponential decay:
  The solution $(s,u)$ exists globally in time, with $ \lim_{t\to\infty}s(t)<\infty$, 
   and $\|u(t)\|_{L^\infty(0,s(t))}$ decays exponentially as $t\to\infty$;
  \item[(b)]  
  Global-in-time solutions with non-exponential decay:
  The solution $(s,u)$ exists globally in time, with $\lim_{t\to\infty}s(t)=\infty$,  
  and $\|u(t)\|_{L^\infty(0,s(t))}$ decays as $t\to\infty$ but not exponentially;
  \item[(c)] Finite-time blow-up solutions:
  $T_m < \infty$ and \eqref{eq:1.2} holds.
\end{itemize}
The occurrence of each type is determined by the size of the initial function $\varphi$.
For any fixed $s_0\in(0,\infty)$, a small initial function $\varphi$ gives rise to type~(1), 
an intermediate one to type~(2), and a large one to type~(3).

This paper is motivated by the papers~\cites{FS, GST, S}, 
which investigated a one-dimensional one-phase Stefan problem for the semilinear heat equation:
\[
\tag{SH}
\label{SH}
\left\{
\begin{array}{ll}
\partial_t v=\partial_x^2 v+v^q, & (t,x)\in D_\sigma(T),\vspace{3pt}\\
-\partial_xv(t,0)=0, &t\in (0,T], \vspace{3pt}\\
v(t,\sigma(t))=0, & t\in (0,T],\vspace{3pt}\\
\sigma'(t)=-\partial_xv(t,\sigma(t)),\qquad & t\in (0,T],\vspace{3pt}\\
v(0,x)=\varphi(x), & x\in [0,\sigma_0],\vspace{3pt}\\
\sigma(0)=\sigma_0,
\end{array}
\right.
\]
where $T\in(0,\infty]$, $q>1$, and $(\sigma_0,\varphi)\in W$ with $\varphi\in C^1([0,s_0])$ and $\varphi'(0)=0$. 
It was shown in \cites{FS, GST, S} that 
all solutions $(\sigma,v)$ to problem~\eqref{SH} fall into one of the following three distinct types:
\begin{itemize}
    \item[(a')] Global-in-time solutions with exponential decay: 
    The solution $(\sigma,v)$ exists globally in time, with 
    $ \lim_{t\to\infty}\sigma(t)<\infty$, 
    and $\|v(t)\|_{L^\infty(0,\sigma(t))}$ decays exponentially as $t\to\infty$;
    \item[(b')] Global-in-time solutions with non-exponential decay: 
    The solution $(\sigma,v)$ exists globally in time, with 
    $\lim_{t\to\infty}\sigma(t)=\infty$,  
    and $\|v(t)\|_{L^\infty(0,\sigma(t))}$ decays as $t\to\infty$ but not exponentially. 
    In particular, 
    $$
    \sigma(t)=O\left(t^{\frac{2}{3}}\right)\mbox{ as }t\to\infty,\qquad
    \liminf_{t\to\infty}\sigma(t)^{\frac{2}{p-1}}\|v(t)\|_{L^\infty(0,\sigma(t))}>0,
    $$
    and hence
    $$
    \liminf_{t\to\infty}t^{\frac{4}{3(p-1)}}\|v(t)\|_{L^\infty(0,\sigma(t))}>0;
    $$
    \item[(c')] Finite-time blow-up solutions: 
    There exists $T\in(0,\infty)$ such that 
    the solution~$(\sigma,v)$ exists on $(0,T)$ and 
    $$
    \limsup_{t\nearrow T}\|v(t)\|_{L^\infty(0,\sigma(t))}=\infty.
    $$ 
\end{itemize}
Similar results to those in \cites{FS, GST, S} have been obtained for one-phase and two-phase Stefan problems for the heat equation with various nonlinear terms; 
see, for example, \cites{BDK, DL, DMZ, LCY, Nin, LW, WC, WZ, ZCX, ZL, ZBL, ZZ} and the references therein. 
We also refer the reader to \cites{A, A2, AI2, AI} for related results on problem~\eqref{SH}.

The main purpose of this paper is to establish a similar trichotomy for problem~\eqref{SP} and 
to clarify how the nonlinear boundary condition affects the transition 
between global-in-time solutions with exponential decay, global-in-time solutions with non-exponential decay, and finite-time blow-up solutions. 
To the best of our knowledge, this is the first paper to address the large-time behavior of solutions to problem~\eqref{SP}; 
for related results without the Stefan boundary condition, see, for example, \cites{FFM, IK} and Remark~\ref{Remark:1.2}. 
\vspace{5pt}

Define 
\begin{equation}
\label{eq:1.4}
\tilde{W}:=\{(a,\phi)\,:\,a\in(0,\infty),\,\,\,\phi\in W^{1,2}(0,a)\setminus\{0\}\mbox{ with $\phi\ge 0$ on $[0,a]$ and $\phi(a)=0$}\}.
\end{equation}
Observe that for each $a\in(0,\infty)$, the Sobolev embedding theorem yields 
\begin{equation}
\label{eq:1.5}
W^{1,2}(0,a)\subset C^\alpha([0,a])
\end{equation} 
for any $\alpha\in(0,1/2)$. 
We give the definition of a solution to problem~\eqref{SP} with initial data $(s_0,\varphi)\in \tilde{W}$. 
%%%
\begin{definition}
\label{Definition:1.1}
Let $(s_0,\varphi)\in \tilde{W}$ and $T\in(0,\infty)$.
We say that a pair $(s,u)$ is a solution to problem~\eqref{SP} on $(0,T]$ if
$s=s(t)\in C([0,T])\cap C^1((0,T])$ is positive on $[0,T]$, 
$$
u\in C^{1;2}(D_s(T)) \cap\, C^{0;1}(D_s^*(T)) \cap\, C(D_s^{**}(T)),
\qquad 
u\ge 0\quad\mbox{in}\quad C(D_s^{**}(T)),
$$
and if $(s,u)$ satisfies all the relations in \eqref{SP} pointwise. 
Here, $D_s(T)$ is defined as in \eqref{eq:1.1}, and 
$$
D_s^*(T):=\bigcup_{t\in(0,T]}\{t\}\times[0,s(t)],\quad
D_s^{**}(T):=\bigcup_{t\in[0,T]}\{t\}\times[0,s(t)].
$$
For any $T\in(0,\infty]$, 
we say that $(s,u)$ is a solution to problem~\eqref{SP} on $(0,T)$ 
if $(s,u)$ is a solution to problem~\eqref{SP} on $(0,T']$ for all $T'\in(0,T)$.  
\end{definition}
%%%
Existence, uniqueness, and the comparison principle for problem~\eqref{SP} are established in Section~2.
%\vspace{3pt}

We are now in a position to present the main results of this paper.
Theorem~\ref{Theorem:1.1} shows that the energy of global-in-time solutions is positive for $t>0$, 
and it also describes the blow-up behavior of solutions exhibiting finite-time blow-up.
\begin{theorem}
\label{Theorem:1.1}
Let $(s,u)$ be a solution to problem~\eqref{SP} with initial data $(s_0,\varphi)\in W$.
Denote by~$T_m$ the maximal existence time of the solution. 
Define 
\begin{equation}
\label{eq:1.6}
E(s(t),u(t)):=\frac{1}{2}\int_0^{s(t)} |\partial_x u(t)|^2\,\dee x-\frac{1}{p+1}u(t,0)^{p+1}\quad\mbox{for}\quad t\in(0,T_m).
\end{equation}
\begin{itemize}
  \item[{\rm (1)}] 
  If $T_m=\infty$, then 
  $$
  E(s(t),u(t))\ge\frac{\pi^2}{256}\frac{\|u(t)\|_{L^1(0,s(t))}^3}{(s(t)+\|u(t)\|_{L^1(0,s(t))})^4}>0\quad\mbox{for $t\in(0,\infty)$}. 
  $$
  \item[{\rm (2)}] 
  If $T_m<\infty$, then $\lim_{t\nearrow T_m}u(t,0)=\infty$ and $\lim_{t\nearrow T_m}s(t)<\infty$. 
  Furthermore, 
  \begin{align}
  \label{eq:1.7}
   & \sup_{t\in(T_m/2,T_m)}\,(T_m-t)^{\frac{1}{2(p-1)}}u(t,0)<\infty,\\
  \label{eq:1.8}
   & \sup_{t\in(T_m/2,T_m)}\|u(t)\|_{L^\infty(\delta,s(t))}<\infty\quad\mbox{for any $\delta\in(0,s_0)$},\\
 \label{eq:1.9}
   & \lim_{t\nearrow T_m}E(s(t),u(t))=-\infty.
  \end{align}
\end{itemize}
\end{theorem}
\begin{remark}
\label{Remark:1.1}
{\rm (1)} 
Even if the initial function $\varphi$ is smooth on $[0,s_0]$ and satisfies the compatibility condition at $x=0$, 
namely $-\partial_x\varphi(0)=\varphi(0)^p$, 
the rescaled function $\lambda\varphi$ with $\lambda>0$ does not necessarily satisfy the compatibility condition. 
Therefore, in our analysis, we do not impose the compatibility condition at $x=0$ on the initial function $\varphi$.
\vspace{2pt}
\newline
{\rm (2)} 
Let $(s,u)$ be a solution to problem~\eqref{SP} with initial data $(s_0,\varphi)\in W$. 
Even if $\varphi\in C^1([0,s_0])$,
in the absence of the compatibility condition,
one cannot expect $\partial_x u$ to be continuous at $(t,x)=(0,0)$.
Therefore, Theorem~{\rm\ref{Theorem:1.1}~(1)} does not necessarily imply that the solution blows up in finite time, 
even when $\varphi\in C^1([0,s_0])$ and $E(s_0,\varphi)<0$.
\end{remark}

Theorem~\ref{Theorem:1.2}, which is the main result of this paper, 
provides a classification of the behavior of solutions to problem~\eqref{SH}.
% into the following three distinct types: 
%global-in-time solutions with exponential decay, global-in-time solutions with non-exponential decay, and finite-time blow-up solutions, depending on the size of the initial function.
%
\begin{theorem}
\label{Theorem:1.2}
Let $(s_0,\varphi)\in W$.
For any $\lambda>0$, let $(s_\lambda,u_\lambda)$ be the solution to problem~\eqref{SP} with initial data $(s_0,\lambda\varphi)$, 
and let $T_\lambda$ denote its maximal existence time.
Then there exist $\lambda_*$, $\lambda^*\in(0,\infty)$ with $\lambda_*\le \lambda^*$ 
such that the following properties hold:
\begin{itemize}
\item[{\rm (1)}] {\bf Global existence with exponential decay}:\vspace{3pt}\\ 
If $0<\lambda<\lambda_*$, then $T_\lambda=\infty$, $\lim_{t\to\infty}s_\lambda(t)<\infty$, and 
$$
\|u_\lambda(t)\|_{L^\infty(0,s_\lambda(t))}=O(e^{-\alpha t})\quad\mbox{as $t\to\infty$ for some $\alpha\in(0,\infty)$}.
$$
\item[{\rm (2)}] {\bf Global existence with non-exponential decay}:\vspace{3pt}\\ 
If $\lambda_*\le\lambda\le \lambda^*$, then $T_\lambda=\infty$ and $\lim_{t\to\infty}s_\lambda(t)=\infty$. 
Furthermore, 
\begin{align*}
 & \lim_{t\to\infty}\|u_\lambda(t)\|_{L^\infty(0,s_\lambda(t))}=0,\quad 
\liminf_{t\to\infty}\,s_\lambda(t)^{\frac{1}{p-1}}\|u_\lambda(t)\|_{L^\infty(0,s_\lambda(t))}>0,\\
 & s_\lambda(t)=(1+o(1))\int_0^t u_\lambda(\tau,0)^p\,\dee\tau=o\left(t^{\frac{1}{2}}\right)\quad\mbox{as $t\to\infty$}.
\end{align*}
In particular, 
\begin{equation}
\label{eq:1.10}
\liminf_{t\to\infty}t^{\frac{1}{2(p-1)}}\|u_\lambda(t)\|_{L^\infty((0,s(t))}=\infty.
\end{equation}
\item[\rm{(3)}] {\bf Blow-up}:\vspace{3pt}\\
If $\lambda>\lambda^*$, the solution~$(s_\lambda,u_\lambda)$ blows up in finite time, i.e., $T_\lambda<\infty$.
\end{itemize}
\end{theorem}
\begin{remark}
\label{Remark:1.2}
{\rm (1)} Let us consider the problem
\begin{equation}
\tag{NBC}
\label{eq:NBC}
\left\{
\begin{array}{ll}
\partial_t v=\partial_x^2 v,\qquad & t>0,\,\,\,x\in(0,\infty),\vspace{3pt}\\
-\partial_x v(t,0)=v(t,0)^p,\qquad & t>0,\vspace{3pt}\\
v(0,x)=\lambda\varphi(x)\ge 0,\qquad & x\in(0,\infty),
\end{array}
\right.
\end{equation}
where $p>1$, $\lambda>0$, and $\varphi\in L^\infty(0,\infty)\cap L^2((0,\infty),e^{|x|^2/4}\,\dee x)\setminus\{0\}$. 
In the case $1<p\le 2$, problem~\eqref{eq:NBC} admits no global-in-time solutions {\rm({\it see, e.g.,} \cites{DFL, GL})}. 
In the case $p>2$, there exists $\lambda_*'>0$ with the following properties  {\rm({\it see, e.g.,} \cite{IK})}.
\begin{itemize}
  \item[{\rm (a)}] If $0<\lambda<\lambda_*'$, then problem~\eqref{eq:NBC} admits a global-in-time solution~$v_\lambda$ such that 
 \begin{equation}
 \label{eq:1.11}
  0<\liminf_{t\to\infty}t^{\frac{1}{2}}\|v_\lambda(t)\|_{L^\infty(0,\infty)}\le \limsup_{t\to\infty}t^{\frac{1}{2}}\|v_\lambda(t)\|_{L^\infty(0,\infty)}<\infty;
  \end{equation}
  \item[{\rm (b)}] If $\lambda=\lambda_*'$, then problem~\eqref{eq:NBC} admits a global-in-time solution~$v_\lambda$ such that
 \begin{equation}
 \label{eq:1.12}
  0<\liminf_{t\to\infty}t^{\frac{1}{2(p-1)}}\|v_\lambda(t)\|_{L^\infty(0,\infty)}\le \limsup_{t\to\infty}t^{\frac{1}{2(p-1)}}\|v_\lambda(t)\|_{L^\infty(0,\infty)}<\infty;
 \end{equation}
  \item[{\rm (c)}] If $\lambda>\lambda_*'$, then the solution~$v_\lambda$ blows up in finite time. 
\end{itemize}
In comparison with the results in Theorem~{\rm\ref{Theorem:1.2}}, 
we see that the decay rates of global-in-time solutions to problem~\eqref{eq:NBC} differ completely from those to problem~\eqref{SP}.
\vspace{3pt}
\newline
{\rm (2)} 
Let $(s_0,\varphi)\in W$, and let $\lambda_*$, $\lambda^*$, and $(s_\lambda,u_\lambda)$ be as in Theorem~{\rm\ref{Theorem:1.2}}.
For any $\lambda>0$, the comparison principle implies that $u_\lambda\le v_\lambda$ in $D_{s_\lambda}(T_{v_\lambda})$, 
where $T_{v_\lambda}$ denotes the maximal existence time of the solution $v_\lambda$ to problem~\eqref{eq:NBC}.
Then we observe from \eqref{eq:1.10} that, 
if $\lambda\ge\lambda_*$, then the solution $v_\lambda$ to problem~\eqref{eq:NBC} does not satisfy both of \eqref{eq:1.11} and \eqref{eq:1.12}, 
that is, the solution $v_\lambda$ blows up in finite time. Hence, we conclude that $\lambda_*>\lambda_*'$.
\end{remark}

Next, under the assumption that $\partial_x\varphi\le 0$ on $[0,s_0]$, 
we obtain a lower bound on the growth rate of $s_\lambda(t)$ as $t\to\infty$ for $\lambda\in[\lambda_*,\lambda^*]$. 
\begin{theorem}
\label{Theorem:1.3}
Let $(s_0,\varphi)\in W$ with $\varphi\in C^1(0,s_0)$ and $\partial_x\varphi\le 0$ on $(0,s_0)$.
Then, for any $\lambda\in[\lambda_*,\lambda^*]$,  
\begin{equation}
\label{eq:1.13}
\liminf_{t\to\infty}t^{-\frac{p-1}{2p-1}}s_{\lambda}(t)>0.
\end{equation}
\end{theorem}
We conjecture that 
\eqref{eq:1.13} holds even without the assumption that $\partial_x\varphi\le 0$ on $(0,s_0)$
and that 
$$
\limsup_{t\to\infty}t^{-\frac{p-1}{2p-1}}s_{\lambda}(t)<\infty,
$$
but both questions remain open.
\vspace{3pt}

We employ techniques from the blow-up analysis for the heat equation with a nonlinear boundary condition 
to establish Theorem~\ref{Theorem:1.1}. (See Propositions~\ref{Proposition:3.1}, \ref{Proposition:3.3}, and \ref{Proposition:4.3}.)
In the proof of Theorem~\ref{Theorem:1.2}, for any $(s_0,\varphi)\in W$, we show that 
$$
\Lambda^*:=\{\lambda>0\,:\,T_\lambda=\infty\}
$$
is bounded, where $T_\lambda$ is as in Theorem~\ref{Theorem:1.2}. 
As in related works, 
we establish this by exploiting the fact that if the energy of the initial data is negative, 
then the corresponding solution blows up in finite time.
However, as stated in Remark~\ref{Remark:1.1}, this strategy does not apply directly to our problem~\eqref{SP}, even when $\varphi\in C^1([0,s_0])$.
Instead, when $\varphi(0)>0$, we construct a suitable initial function that satisfies the compatibility condition and whose solution blows up in finite time, thereby proving the boundedness of $\Lambda^*$.
When $\varphi(0)=0$, by considering the heat equation with the homogeneous Neumann boundary condition, we reduce the problem to the case $\varphi(0)>0$ and again obtain the boundedness of $\Lambda^*$.
The remaining assertions of Theorem~\ref{Theorem:1.2} are obtained by suitable modifications of the arguments in \cites{FS, GST, S}, together with 
a result on the growth rate of the free boundary in the one-phase Stefan problem for the heat equation with an inhomogeneous Dirichlet boundary condition.
We note that several parts of the proofs require a delicate analysis to handle the nonlinear boundary condition 
(see, for example, the proof of Propositions~\ref{Proposition:2.2} and \ref{Proposition:4.1}). 
For the proof of Theorem~\ref{Theorem:1.3},
under the assumption that $\partial_x\varphi\le 0$ on $(0,s_0)$, 
we use the relations between the solutions $(s_\lambda,u_\lambda)$ in Theorem~\ref{Theorem:1.2} to derive an inequality for the integral $\int_0^t u_\lambda(\tau,0)^p\,\dee\tau$.
This, in turn, yields lower bounds on the growth rate of $s_\lambda(t)$ as $t\to\infty$, thereby completing the proof of Theorem~\ref{Theorem:1.3}.

The rest of this paper is organized as follows.
In Section~2, we present preliminary results on problem~\eqref{SP}, including existence, uniqueness, stability, the comparison principle, and energy estimates.
Section~3 provides a sufficient condition for finite-time blow-up and describes the behavior of blow-up solutions.
In Section~4, we establish several propositions on global-in-time solutions; in particular, we show that for any global-in-time solution $(s,u)$,
the $L^\infty$-norm of $u$ converges to zero as $t\to\infty$, and we give a sufficient condition for global existence with exponential decay.
We also prove that, if the solution blows up in finite time, then its energy diverges at the blow-up time.
Finally, in Section~5, we complete the proofs of Theorems~\ref{Theorem:1.1} and~\ref{Theorem:1.2}. 
Furthermore, we prove Theorem~\ref{Theorem:1.3}. 
%%%%%%%%%%%%%%%%%%%%%%%%%%%%%%
%%%%%%%%%%%%%%%%%%%%%%%%%%%%%%
\section{Preliminaries}
%%%%%%%%%%%%%%%%%%%%%%%%%%%%%%
%%%%%%%%%%%%%%%%%%%%%%%%%%%%%%
In this section, we collect some preliminary results on solutions to problem~\eqref{SP}.
Throughout this paper, for any function $f$ defined on an interval $I \subset \mathbb{R}$, 
we extend $f$ by zero outside $I$ and, when convenient, identify $f$ with its zero extension.
For simplicity, we also write
$$
\|f\|_{L^r}:=\|f\|_{L^r(I)}
$$
for $r\in[1,\infty]$ whenever no confusion arises.
We use $C$ to denote generic positive constants, which may vary from line to line.
Furthermore, for any Lipschitz continuous function $g$ on an interval $I\subset{\mathbb R}$, we define 
$$
\|g\|_{{\rm Lip}(I)}:=\sup_{x\in I}|g(x)|+\sup_{\substack{x,y\in I \\ x\not=y}}\frac{|g(x)-g(y)|}{|x-y|}.
$$
We write $D_s:=D_s(\infty)$ for simplicity. 

We first establish the local-in-time existence of solutions to problem~\eqref{SP} with initial data in~$\tilde{W}$, where $\tilde{W}$ is defined in \eqref{eq:1.4}.
\begin{proposition}
\label{Proposition:2.1} 
Let $M>0$. 
Then there exists $T\in(0,\infty)$ such that, for any $(s_0,\varphi)\in \tilde{W}$ satisfying 
\begin{equation}
\label{eq:2.1}
s_0+s_0^{-1}+\|\varphi\|_{W^{1,2}([0,s_0])}\le M, 
\end{equation}
problem~\eqref{SP} with initial data $(s_0,\varphi)$ admits a solution $(s,u)$ on $(0,T]$ with
\begin{equation}
\label{eq:2.2}
\sup_{t\in(0,T]}\|u(t)\|_{L^\infty}+\sup_{t\in(0,T]}s(t)\le 2(M^{\frac{3}{2}}+M+1). 
\end{equation}
\end{proposition}
{\bf Proof.}
Let $M>0$, and assume \eqref{eq:2.1}. Set $M':=M^{3/2}+1$.
Let $F$ be a non-decreasing, $C^1$-smooth function on $[0,\infty)$ such that 
\begin{equation*}
\begin{array}{ll}
F(\eta)=\eta^p\quad\mbox{for $\eta\in[0,2M')$},
\quad
 & F(\eta)\le\eta^p\quad\mbox{for}\quad \eta\in[2M',2M'+2],\vspace{5pt}\\
F(\eta)=(2M'+1)^p\quad\mbox{for}\quad \eta\in(2M'+2,\infty). & 
\end{array}
\end{equation*}
By \cite{K}*{Theorem~2.2} (see also \cite{FP}) and \eqref{eq:1.5}, 
there exists $T_1\in(0,\infty)$, depending only on $F$, $M$, and $p$, 
such that the one-phase Stefan problem with nonlinear boundary condition
$$
\left\{
\begin{array}{ll}
\partial_t u=\partial_x^2 u, & (t,x)\in D_s(T_1),\vspace{3pt}\\
-\partial_xu(t,0)=F(|u(t,0)|),\qquad & t\in(0,T_1], \vspace{3pt}\\
u(t,s(t))=0, &  t\in(0,T_1],\vspace{3pt}\\
s'(t)=-\partial_xu(t,s(t)), &  t\in(0,T_1], \vspace{3pt}\\
u(0,x)=\varphi(x), & x\in [0,s_0],\vspace{3pt}\\
s(0)=s_0,
\end{array}
\right.
$$
admits a classical solution $(s,u)$ on $(0,T_1]$ such that 
\begin{equation}
\label{eq:2.3}
s_0<s(t)<2s_0\le 2M,\quad t\in[0,T_1].
\end{equation}
Furthermore, by the maximum principle and Hopf's lemma, we obtain $u\ge 0$ in $D_s^{**}(T_1)$. 
On the other hand, since 
$$
0\le\varphi(x)=-\int^{s_0}_x\partial_x\varphi\,\dee x\le s_0^{\frac{1}{2}}\|\partial_x\varphi\|_{L^2}\le M^{\frac{3}{2}},\quad x\in[0,s_0],
$$
extending $\varphi$ by zero to $(0,\infty)$, 
we find $T_2\in(0,\infty)$, depending only on $M$ and $p$, 
such that the problem
\begin{equation}
\label{eq:2.4}
\left\{
\begin{array}{ll}
\partial_t v=\partial_x^2 v, & (t,x)\in (0,T_2]\times(0,\infty),\vspace{3pt}\\
-\partial_xv(t,0)=v(t,0)^p, & t\in(0,T_2],\vspace{3pt}\\
v(0,x)=\varphi(x)+1, & x\in [0,\infty),
\end{array}
\right.
\end{equation}
admits a positive classical solution $v$ satisfying
$$
0<v(t,x)\le 2(M^{\frac{3}{2}}+1)=2M',\quad (t,x)\in(0,T_2]\times(0,\infty).
$$
(See, e.g., \cites{IS, IS2} for the existence of bounded solutions to problem~\eqref{eq:2.4}, 
and  \cite{LSU}*{Chapter~IV} together with \cite{DiB}*{Theorem~6.2} and \cite{IK}*{Lemma~2.4} for regularity results.)
Set
$$
w(t,x):=v(t,x)-u(t,x),\quad (t,x)\in D_s(T_3),
\quad\mbox{where}\quad T_3=\min\{T_1,T_2\}.
$$
Since 
$$
-\partial_x w(t,0)=v(t,0)^p-F(|u(t,0)|)\ge v(t,0)^p-u(t,0)^p=a(t)w(t,0),\quad t\in(0,T_3],
$$
where 
\begin{equation}
\label{eq:2.5}
a(t):=
\left\{
\begin{array}{ll} 
\displaystyle{\frac{v(t,0)^p-u(t,0)^p}{v(t,0)-u(t,0)}} & \quad\mbox{if}\quad v(t,0)\not=u(t,0),\vspace{5pt}\\
pv(t,0)^{p-1} & \quad\mbox{if}\quad v(t,0)=u(t,0),
\end{array}
\right.
\end{equation}
$w$ satisfies
\begin{equation*}
\left\{
\begin{array}{ll}
\partial_t w=\partial_x^2 w, & (t,x)\in D_s(T_3),\vspace{3pt}\\
-\partial_xw(t,0)\ge a(t)w(t,0), & t\in(0,T_3],\vspace{3pt}\\
w(t,s(t))>0, & t\in(0,T_3],\vspace{3pt}\\
w(0,x)=1, & x\in [0,s_0].
\end{array}
\right.
\end{equation*}
The maximum principle, together with Hopf's lemma, yields 
\begin{equation}
\label{eq:2.6}
w>0\quad\mbox{in}\quad D_s^{**}(T_3). 
\end{equation}
Indeed, if $w$ is nonpositive for some point $(t_0,x_0)\in D_s^{**}(T_3)$, then, 
by the continuity of $w$, we can find $(t_1,x_1)\in D_s^*(T_3)$ such that  
$$
w(t_1,x_1)=0,\qquad \mbox{$w(t,x)>0$ for $(t,x)\in D_s^{**}(T_3)$ with $t<t_1$}.
$$
It follows from the maximum principle that $x_1=0$ and $w>0$ in $D_s^{**}\setminus\{(t_1,0)\}$, which, together with Hopf's lemma, implies that 
$0>-\partial_x w(t_1,0)\ge a(t_1)w(t_1,0)=0$. 
This is a contradiction, and hence \eqref{eq:2.6} holds. 
Consequently,
\begin{equation}
\label{eq:2.7}
0\le u(t,x)<v(t,x)\le 2M',\quad (t,x)\in D_s^{**}(T_3).
\end{equation}
Then $(s,u)$ is a solution to problem~\eqref{SP} on $(0,T_3]$. 
Furthermore, by \eqref{eq:2.3} and \eqref{eq:2.7}, 
we obtain \eqref{eq:2.2}. Thus, Proposition~\ref{Proposition:2.1} follows.
$\Box$\vspace{5pt}

We obtain inequalities on the energy $E(s(t),u(t))$ and the $L^1$-norm $\|u(t)\|_{L^1}$. 
Lemma~\ref{Lemma:2.1} plays a crucial role in the proofs of Theorem~\ref{Theorem:1.1}~(1) and Theorem~\ref{Theorem:1.2}.
\begin{lemma}
\label{Lemma:2.1}
Let $(s,u)$ be a solution to problem~\eqref{SP} on $(0,T)$ with initial data $(s_0,\varphi)\in W$, where $T\in(0,\infty]$. 
Let $E(s(t),u(t))$ be defined as in \eqref{eq:1.6} for $t\in(0,T)$.
Then
\begin{align}
\label{eq:2.8}
 & \frac{\dee}{\dee t}E(s(t),u(t))=-\int_0^{s(t)} |\partial_t u(t)|^2\,\dee x-\frac{1}{2}s'(t)^3,\\
\label{eq:2.9}
 & \|u(t)\|_{L^1}-\|\varphi\|_{L^1}=s_0-s(t)+\int_0^tu(\tau,0)^p\,\dee\tau,
\end{align}
for $t\in(0,T)$. In particular, 
\begin{equation}
\label{eq:2.10}
\int_{t_1}^{t_2}\int_0^{s(t)} |\partial_t u|^2\,\dee x\,\dee t+\frac{1}{2}\int_{t_1}^{t_2} (s')^3\,\dee t
=E(s(t_1),u(t_1))-E(s(t_2),u(s_2))
\end{equation}
for $t_1$, $t_2\in(0,T)$ with $t_1\le t_2$.
\end{lemma}
{\bf Proof.}
Let $(s,u)$ be a solution to problem~\eqref{SP} on $(0,T)$ with initial data $(s_0,\varphi)\in W$, where $T\in(0,\infty]$.
For any $t\in(0,T)$, 
since $u(t,s(t))=0$ and $\partial_xu(t,s(t))=-s'(t)$, 
we have
$$
0=\frac{\dee}{\dee t}u(t,s(t))=\partial_t u(t,s(t))+\partial_x u(t,s(t))s'(t)=\partial_t u(t,s(t))-s'(t)^2,
$$
which implies that
\begin{equation}
\label{eq:2.11}
\partial_x u(t,s(t)) \partial_tu(t,s(t))=-s'(t)^3,\quad t\in(0,T).
\end{equation}
Then 
\begin{align*}
\frac{\dee}{\dee t}E(s(t),u(t)) &=\int_0^{s(t)}\partial_xu(t)\partial_t\partial_xu(t)\,\dee x
+\frac{1}{2}|\partial_x u(t,s(t))|^2 s'(t)-u(t,0)^p \partial_tu(t,0)\\
&=\partial_xu(t,s(t))\partial_t u(t,s(t))-\partial_x u(t,0)\partial_tu(t,0)\\
& \qquad\qquad
-\int_0^{s(t)}|\partial_t u(t)|^2\,\dee x
+\frac{1}{2}s'(t)^3-u(t,0)^p \partial_tu(t,0)\\
&=-\int_0^{s(t)}|\partial_t u(t)|^2\,\dee x-\frac{1}{2}s'(t)^3
\end{align*}
for $t\in(0,T)$, 
which is exactly~\eqref{eq:2.8}. Then we immediately obtain \eqref{eq:2.10}. 
Furthermore, since
\begin{align*}
\frac{\dee}{\dee t}\int_0^{s(t)}u(t)\,\dee x&=\int_0^{s(t)}\partial_tu(t)\,\dee x+u(t,s(t))s'(t)=\int_0^{s(t)} \partial_x^2u(t)\,\dee x\\
&=\partial_xu(t,s(t))-\partial_xu(t,0)=-s'(t)+u(t,0)^p,
\end{align*}
we have
$$
\int_0^{s(t)}u(t)\,\dee x-\int_0^{s_0}u(0)\,\dee x=s_0-s(t)+\int_0^t u(\tau,0)^p\,\dee \tau
$$
for $t\in(0,T)$. 
This implies~\eqref{eq:2.9}, and thus Lemma~\ref{Lemma:2.1} follows.
$\Box$\vspace{5pt}

We next obtain pointwise estimates of solutions to problem~\eqref{SP}.
\begin{lemma}
\label{Lemma:2.2}
Let $(s,u)$ be a solution to problem~\eqref{SP} on $(0,T]$ with initial data $(s_0,\varphi)\in W$, 
where $T\in(0,\infty)$. 
Assume that there exists $M>0$ satisfying 
\begin{equation}
\label{eq:2.12}
\begin{array}{ll}
0\le\varphi(x)\le M(s_0-x), & \quad x\in[0,s_0],\vspace{5pt}\\
0\le u(t,0)^p\le M, & \quad t\in(0,T].
\end{array}
\end{equation}
Then
$$
0\le u(t,x)\le M(s(t)-x),\quad (t,x)\in D_s^{**}(T).
$$
In particular, 
\begin{equation}
\label{eq:2.13}
s'(t)=-\partial_xu(t,s(t))\le M,\qquad t\in(0,T].
\end{equation}
\end{lemma}
{\bf Proof.}
Set 
$$
v(t,x):=M(s(t)-x),\quad (t,x)\in D_s(T).
$$
By \eqref{eq:1.3}, $v$ satisfies 
$$
\left\{
\begin{array}{ll}
\partial_t v=\partial_x^2 v+Ms'(t)\ge \partial_x^2 v, & (t,x)\in D_s(T),\vspace{3pt}\\
-\partial_xv(t,0)=M,\qquad & t\in(0,T], \vspace{3pt}\\
v(t,s(t))=0, &  t\in(0,T],\vspace{3pt}\\
v(0,x)=M(s_0-x), & x\in [0,s_0],\vspace{3pt}\\
s(0)=s_0.
\end{array}
\right.
$$
The comparison principle, together with \eqref{eq:2.12}, implies that 
$$
0\le u(t,x)\le v(t,x)=M(s(t)-x),\quad (t,x)\in D_s^{**}(T). 
$$
Then \eqref{eq:2.13} also holds, and the proof of Lemma~\ref{Lemma:2.2} is complete. 
$\Box$\vspace{5pt}

We employ Lemma~\ref{Lemma:2.2} to obtain a result on the stability of solutions to problem~\eqref{SP}.
\begin{proposition}
\label{Proposition:2.2}
Let $(s,u)$ and $(\sigma,v)$ be solutions to problem~\eqref{SP} on $(0,T]$ with initial data $(s_0,\varphi)\in W$ and $(\sigma_0,\psi)\in W$, respectively, 
where $T\in(0,\infty)$. 
Let $M>0$ be such that 
\begin{align}
\label{eq:2.14}
 & 0\le\varphi(x)\le M(s_0-x),\quad x\in [0,s_0],\\
\label{eq:2.15}
  & 0\le\psi(x)\le M(\sigma_0-x),\quad x\in [0,\sigma_0].
\end{align}
Then, for any $\epsilon>0$, 
there exists $\eta>0$ such that 
if 
\begin{equation}
\label{eq:2.16}
|s_0-\sigma_0|+\|\varphi-\psi\|_{L^\infty}\le \eta, 
\end{equation}
then
\begin{equation}
\label{eq:2.17}
\sup_{t\in(0,T]}|s(t)-\sigma(t)|+\sup_{t\in(0,T]}\|u(t)-v(t)\|_{L^\infty}\le \epsilon.
\end{equation}
Here $\eta$ depends only on $\epsilon$, $p$, $s_0$, $\sup_{t\in(0,T]}u(t,0)$, $T$, and $M$. 
\end{proposition}
{\bf Proof.}
Suppose that $M>0$ is chosen so that \eqref{eq:2.14} and \eqref{eq:2.15} hold.
Let $\eta\in(0,\min\{1,s_0/2\})$ be sufficiently small, and assume \eqref{eq:2.16}.
Choose $\delta>0$ satisfying
\begin{equation}
\label{eq:2.18}
\delta<\min\left\{1,\frac{s_0^2}{16}\right\}.
\end{equation}
Let $T_\delta:=\min\{T,\delta\}$.
In the proof, generic constants $C$ are independent of $\eta$ and $\delta$. 
We identify $u(t)$ and $v(t)$ with their zero extensions, respectively. 
Set
$$
w(t,x):=u(t,x)-v(t,x),\quad \alpha(t):=\min\{s(t),\sigma(t)\},\quad \beta(t):=\max\{s(t),\sigma(t)\},
$$
for $t\in(0,T]$ and $x\in[0,\infty)$.  
Define 
$$
T_*:=\sup\left\{t\in(0,T_\delta]\,:\, \mu(\tau)+\rho(\tau)\le 1\mbox{ for $\tau\in[0,t]$}\right\},
$$
where
\begin{equation}
\label{eq:2.19}
\mu(t):=\sup_{\tau\in(0,t]}\|u(\tau)-v(\tau)\|_{L^\infty},\quad \rho(t):=\sup_{\tau\in(0,t]}|s(\tau)-\sigma(\tau)|.
\end{equation}
It follows from Definition~\ref{Definition:1.1} that $\sup_{t\in(0,T]}u(t,0)<\infty$. 
Then, by the definition of $T_*$, taking $M>0$ sufficiently large if necessary, 
we obtain 
\begin{equation}
\label{eq:2.20}
\sup_{t\in(0,T_*]}v(t,0)^p\le \left(\sup_{t\in(0,T_*]}u(t,0)+1\right)^p\le M. 
\end{equation}
By \eqref{eq:2.14}, \eqref{eq:2.15}, and \eqref{eq:2.20}, we apply Lemma~\ref{Lemma:2.2} to obtain 
\begin{equation}
\label{eq:2.21}
\begin{split}
 & 0\le u(t,x)\le M(s(t)-x),\quad (t,x)\in D_s^{**}(T_*),\\
 & 0\le v(t,x)\le M(\sigma(t)-x),\quad (t,x)\in D_\sigma^{**}(T_*),\\
 & 0\le s'(t)\le M,\quad 0\le\sigma'(t)\le M,\quad t\in(0,T_*].
\end{split}
\end{equation}
Since $T_*\le\delta\le 1$ and $0\le\eta\le 1$, we have
\begin{align}
\label{eq:2.22}
 & \sup_{t\in(0,T_*]}\alpha(t)\le\sup_{t\in(0,T_*]}\beta(t)\le s_0+1+M,\\
\label{eq:2.23}
 & \rho(t)\le \eta+Mt,\quad t\in(0,T_*].
\end{align}
%\vspace{5pt}
%\newline
%
\underline{Step 1.} We derive estimates of $\mu(t)$ in $D_\alpha(T_*)$.  
Let $t_*\in(0,T_*]$. 
It follows from \eqref{eq:2.20} and \eqref{eq:2.21} that
$w$ satisfies 
\begin{equation}
\label{eq:2.24}
\left\{
\begin{array}{ll}
\partial_t w=\partial_x^2 w, & (t,x)\in D_\alpha(t_*),\vspace{3pt}\\
-\partial_xw(t,0)\le pM^{\frac{p-1}{p}}w(t,0), & t\in(0,t_*],\vspace{3pt}\\
w(t,\alpha(t))\le M|s(t)-\sigma(t)|\le M\rho(t_*), & t\in(0,t_*],\vspace{3pt}\\
w(0,x)\le\|\varphi-\psi\|_{L^\infty}, & x\in[0,\alpha(0)].
\end{array}
\right.
\end{equation}
Define 
\begin{align*}
z(t,x) & :=L(4\pi (t+1))^{-\frac{1}{2}}\exp\left(-\frac{x^2}{4(t+1)}\right)(\|\varphi-\psi\|_{L^\infty}+\rho(t_*))\\
 & \qquad\quad
 +2pM^{\frac{p-1}{p}}\mu(t_*)\int_0^t (4\pi(t-s))^{-\frac{1}{2}}\exp\left(-\frac{x^2}{4(t-s)}\right)\,\dee s
\end{align*}
for $(t,x)\in[0,t_*]\times[0,\infty)$, where $L>0$.  
By \eqref{eq:2.22} and the fact that $t_*\le \delta\le 1$, $z$ satisfies 
\begin{equation*}
\left\{
\begin{array}{ll}
\partial_t z=\partial_x^2 z,\quad & (t,x)\in(0,t_*]\times (0,\infty),\vspace{3pt}\\
-\partial_x z(t,0)=pM^{\frac{p-1}{p}}\mu(t_*),\quad & t\in(0,t_*],\vspace{3pt}\\
\displaystyle{z(t,\alpha(t))\ge CL\rho(t_*)}, \quad & t\in(0,t_*],\vspace{3pt}\\
z(0,x)\ge CL\|\varphi-\psi\|_{L^\infty},\quad & x\in[0,\alpha(0)].
\end{array}
\right.
\end{equation*}
Consequently, 
taking $L>0$ sufficiently large,  we have
\begin{equation*}
\left\{
\begin{array}{ll}
\partial_t z=\partial_x^2z, & (t,x)\in D_\alpha(t_*),\vspace{7pt}\\
-\partial_xz(t,0)=pM^{\frac{p-1}{p}}\mu(t_*),& t\in(0,t_*],\vspace{7pt}\\
z(t,\alpha(t))>M\rho(t_*), & t\in(0,t_*],\vspace{7pt}\\
z(0,x)>\|\varphi-\psi\|_{L^\infty}, & x\in[0,\alpha(0)].
\end{array}
\right.
\end{equation*}
Then it follows from the maximum principle and Hopf's lemma that 
$$
u(t,x)-v(t,x)=w(t,x)<z(t,x)\le C\|\varphi-\psi\|_{L^\infty}+C\rho(t_*)+C\delta^{\frac{1}{2}}M^{\frac{p-1}{p}}\mu(t_*)
$$
for $(t,x)\in D_\alpha^{**}(t_*)$.
In particular, 
\begin{equation}
\label{eq:2.25}
u(t_*,x)-v(t_*,x)\le C\|\varphi-\psi\|_{L^\infty}+C\rho(t_*)+C\delta^{\frac{1}{2}}M^{\frac{p-1}{p}}\mu(t_*)
\end{equation}
for $x\in[0,\alpha(t_*)]$. 
Similarly, we have 
\begin{equation}
\label{eq:2.26}
v(t_*,x)-u(t_*,x)\le C\|\varphi-\psi\|_{L^\infty}+C\rho(t_*)+C\delta^{\frac{1}{2}}M^{\frac{p-1}{p}}\mu(t_*) 
\end{equation}
for $x\in[0,\alpha(t_*)]$.
Since $t_*\in(0,T_*]$ is arbitrary, 
combining \eqref{eq:2.25} and \eqref{eq:2.26}, 
we obtain 
$$
\|u(t)-v(t)\|_{L^\infty(0,\alpha(t))}\le C\|\varphi-\psi\|_{L^\infty}+C\rho(t)+C\delta^{\frac{1}{2}}M^{\frac{p-1}{p}}\mu(t)
$$
for $t\in(0,T_*]$. 
This, together with \eqref{eq:2.21}, implies that
\[
\mu(t) \le \sup_{\tau\in(0,t]}\|u(\tau)-v(\tau)\|_{L^\infty(0,\alpha(t))}+M\sup_{\tau\in(0,t]}|s(\tau)-\sigma(\tau)|
 \le C\eta+C\delta^{\frac{1}{2}}M^{\frac{p-1}{p}}\mu(t)+C\rho(t)
\]
for $t\in(0,T_*]$.
Taking sufficiently small $\delta>0$ if necessary, we obtain 
\begin{equation}
\label{eq:2.27}
\mu(t)\le C\eta+C\rho(t)
\end{equation}
for $t\in(0,T_*]$.
\vspace{5pt}
\newline
\underline{Step 2.} 
We obtain $L^2$-estimates of $w$. 
Since $\eta<s_0/2$, it follows that 
$$
\alpha(0)\ge s_0-\eta>\frac{s_0}{2}.
$$
This, together with \eqref{eq:2.18}, yields 
$$
\delta^{\frac{1}{2}}\le\frac{s_0}{4}<\frac{\alpha(0)}{2}.
$$
Let $\zeta$ be a smooth function on $[0,\infty)$ such that 
\begin{equation}
\label{eq:2.28}
\zeta=1\mbox{ on $[0,\alpha(0)-2\delta^{\frac{1}{2}}]$},\quad \zeta=0\mbox{ on $[\alpha(0)-\delta^{\frac{1}{2}},\infty)$},
\quad |\partial_x\zeta|\le 2\delta^{-\frac{1}{2}}\mbox{ on $[0,\infty)$}.
\end{equation}
Since $\alpha(t)\ge\alpha(0)$ for $t\in[0,T_*]$, 
it follows from \eqref{eq:2.28} that
\begin{align*}
 & \int_0^t\int_0^{\alpha(\tau)} \partial_tw\cdot w\zeta^2\,\dee x\,\dee\tau\\
 & =\frac{1}{2}\int_0^t \left(\frac{\dee}{\dee \tau}\int_0^{\alpha(\tau)} w(\tau)^2\zeta^2\,\dee x\right)\,\dee\tau
=\frac{1}{2}\int_0^{\alpha(t)} w(t)^2\zeta^2\,\dee x-\frac{1}{2}\int_0^{\alpha(0)} w(0)^2\zeta^2\,\dee x
\end{align*}
for $t\in(0,T_*]$.
Furthermore, since 
\begin{equation}
\label{eq:2.29}
-\partial_x w(t,0)w(t,0)\le |u(t,0)^p-v(t,0)^p| |w(t,0)|\le Cw(t,0)^2,\quad t\in(0,T_*],
\end{equation}
by \eqref{eq:2.24} and \eqref{eq:2.28}, we obtain
\begin{align*}
 & \int_0^t\int_0^{\alpha(\tau)} \partial_tw\cdot w\zeta^2\,\dee x\,\dee\tau=\int_0^t\int_0^{\alpha(\tau)} \partial_x^2w\cdot w\zeta^2\,\dee x\,\dee\tau\\
 & =-\int_0^t \partial_x w(\tau,0)w(\tau,0)\,\dee\tau-\int_0^t\int_0^{\alpha(\tau)}|\partial_x w|^2\zeta^2\,\dee x\,\dee\tau
 -2\int_0^t\int_0^{\alpha(\tau)}\partial_x w w\zeta\partial_x\zeta\,\dee x\,\dee\tau\\
 & \le C\int_0^t w(\tau,0)^2\,\dee\tau
 -\frac{1}{2}\int_0^t\int_0^{\alpha(\tau)}|\partial_x w|^2\zeta^2\,\dee x\,\dee\tau
 +C\int_0^t\int_0^{\alpha(\tau)}w^2|\partial_x\zeta|^2\,\dee x\,\dee\tau
\end{align*}
for $t\in(0,T_*]$. These, together with \eqref{eq:2.19}, \eqref{eq:2.22}, \eqref{eq:2.28}, and the fact that $T_*\le\delta\le 1$, imply that 
\begin{align*}
\int_0^{\alpha(0)-2\delta^{\frac{1}{2}}}w(t)^2\,\dee x
 & \le \int_0^{\alpha(0)}w(0)^2\,\dee x+C\int_0^t w(\tau,0)^2\,\dee\tau+\frac{C}{\delta}\int_0^t\int_{\alpha(0)-2\delta^{\frac{1}{2}}}^{\alpha(0)-\delta^{\frac{1}{2}}} w^2\,\dee x\,\dee\tau\\
 & \le C\|\varphi-\psi\|_{L^\infty}^2+C(\delta+\delta^{\frac{1}{2}})\mu(t)^2,\quad t\in(0,T_*].
\end{align*}
Consequently, we obtain
\begin{equation}
\label{eq:2.30}
\left(\int_0^{\alpha(0)-2\delta^{\frac{1}{2}}}w(t)^2\,\dee x\right)^{\frac{1}{2}}
\le C\eta+C\delta^{\frac{1}{4}}\mu(t),\quad t\in(0,T_*].
\end{equation}
\underline{Step 3.} 
We complete the proof. 
It follows from Lemma~\ref{Lemma:2.1} that 
\begin{align*}
 & s(t)=s_0+\int_0^t u(\tau,0)^p\,\dee\tau-\int_0^{s(t)}u(t)\,\dee x+\int_0^{s_0}\varphi\,\dee x,\\
 & \sigma(t)=\sigma_0+\int_0^t v(\tau,0)^p\,\dee\tau-\int_0^{\sigma(t)} v(t)\,\dee x+\int_0^{\sigma_0}\psi\,\dee x,
\end{align*}
for $t\in(0,T_*]$. Then, by \eqref{eq:2.22}, \eqref{eq:2.24}, and \eqref{eq:2.29}, we have
\begin{equation}
\label{eq:2.31}
\begin{split}
|s(t)-\sigma(t)| & \le |s_0-\sigma_0|+C\delta\mu(t)+\int_0^{\alpha(t)}|u(t)-v(t)|\,\dee x\\
 & \quad+\int_{\alpha(t)}^{s(t)}u(t)\,\dee x+\int_{\alpha(t)}^{\sigma(t)}v(t)\,\dee x+\beta(0)\|\varphi-\psi\|_{L^\infty}\\
 & \le |s_0-\sigma_0|+C\delta\mu(t)+C\left(\int_0^{\alpha(0)-2\delta^{\frac{1}{2}}}|w(t)|^2\,\dee x\right)^{\frac{1}{2}}\\
 & \quad +\int_{\alpha(0)-2\delta^{\frac{1}{2}}}^{\alpha(t)}|w(t)|\,\dee x+\int_{\alpha(t)}^{s(t)}u(t)\,\dee x+\int_{\alpha(t)}^{\sigma(t)}v(t)\,\dee x+C\|\varphi-\psi\|_{L^\infty}
\end{split}
\end{equation}
for $t\in(0,T_*]$. 
On the other hand, 
by \eqref{eq:2.21}, if $s(t)\ge\sigma(t)$, then
$$
\int_{\alpha(t)}^{s(t)}u(t)\,\dee x+\int_{\alpha(t)}^{\sigma(t)}v(t)\,\dee x
=\int_{\sigma(t)}^{s(t)}u(t)\,\dee x\le M\int_{\sigma(t)}^{s(t)}(s(t)-x)\,\dee x=\frac{M}{2}(s(t)-\sigma(t))^2,
$$
whereas if $s(t)<\sigma(t)$, then
$$
\int_{\alpha(t)}^{s(t)}u(t)\,\dee x+\int_{\alpha(t)}^{\sigma(t)}v(t)\,\dee x
=\int_{s(t)}^{\sigma(t)}v(t)\,\dee x\le M\int_{s(t)}^{\sigma(t)}(\sigma(t)-x)\,\dee x=\frac{M}{2}(s(t)-\sigma(t))^2.
$$
These, together with \eqref{eq:2.23}, implies that 
\begin{equation}
\label{eq:2.32}
\int_{\alpha(t)}^{s(t)}u(t)\,\dee x+\int_{\alpha(t)}^{\sigma(t)}v(t)\,\dee x\le C\rho(t)^2
\le C(\eta+M\delta)\rho(t)\le C\eta+C\delta\rho(t),\quad t\in(0,T_*].
\end{equation}
Combining \eqref{eq:2.27}, \eqref{eq:2.30}, \eqref{eq:2.31}, and \eqref{eq:2.32}, 
we obtain 
$$
\rho(t)\le C\eta+C(\delta+\delta^{\frac{1}{2}}+\delta^{\frac{1}{4}})\mu(t)+C\delta\rho(t)
\le C\eta+C\delta^{\frac{1}{4}}\rho(t),\quad t\in(0,T_*].
$$
Then, taking $\delta>0$ sufficiently small if necessary, 
we obtain 
$$
\rho(t)\le C\eta,\quad t\in(0,T_*].
$$
Let $\epsilon\in(0,1/2)$. 
Then, by \eqref{eq:2.27}, and taking $\delta>0$ and $\eta>0$ sufficiently small if necessary, 
we obtain 
\begin{equation}
\label{eq:2.33}
\rho(t)+\mu(t)\le C\eta\le \epsilon<\frac{1}{2},\quad t\in(0,T_*],
\end{equation}
which implies that $T_*=T_\delta=\min\{T,\delta\}$. 
Therefore, 
by \eqref{eq:2.33}, 
if $T\le\delta$, then $T_\delta=T$, hence 
$$
\sup_{t\in(0,T]}\|u(t)-v(t)\|_{L^\infty}+\sup_{t\in(0,T]}|s(t)-\sigma(t)|\le \epsilon.
$$
This yields \eqref{eq:2.17}, and Proposition~\ref{Proposition:2.2} follows. 
On the other hand, 
if $T>\delta$, then $T_\delta=\delta$, hence 
$$
\sup_{t\in(0,\delta]}\|u(t)-v(t)\|_{L^\infty}+\sup_{t\in(0,\delta]}|s(t)-\sigma(t)|\le C\eta\le\epsilon. 
$$
Repeating the above argument and taking $\eta>0$ sufficiently small, 
we obtain 
$$
\sup_{t\in(0,T_{2\delta}]}\|u(t)-v(t)\|_{L^\infty}+\sup_{t\in(0,T_{2\delta}]}|s(t)-\sigma(t)|\le C\eta\le\epsilon.
$$
If $T\le 2\delta$, then Proposition~\ref{Proposition:2.2} follows.
By repeating this argument finitely many times, we complete the proof of Proposition~\ref{Proposition:2.2}.
$\Box$\vspace{7pt}

As a corollary of Proposition~\ref{Proposition:2.2}, we deduce the uniqueness of solutions to problem~\eqref{SP}.
\begin{proposition}
\label{Proposition:2.3} 
For any $T\in(0,\infty)$, 
problem~\eqref{SP} with initial data $(s_0,\varphi)\in W$ admits at most one solution $(s,u)$ on $(0,T]$.
\end{proposition}

Next, we show that, if the maximal existence time $T_m$ is finite, 
then \eqref{eq:1.2} holds; in fact,  
\begin{equation}
\label{eq:2.34}
\lim_{t\nearrow T_m}u(t,0)=\infty. 
\end{equation}
\begin{proposition}
\label{Proposition:2.4}
Let $(s,u)$ be a solution to problem~\eqref{SP} with initial data $(s_0,\varphi) \in \tilde{W}$, 
and let $T_m$ denote its maximal existence time.
If $T_m<\infty$, then \eqref{eq:2.34} holds. 
\end{proposition} 
{\bf Proof.}
The proof is by contradiction. 
Assume that \eqref{eq:2.34} does not hold. 
Then there exists $\{t_n\}\subset(0,T_m)$ such that 
\begin{equation}
\label{eq:2.35}
\lim_{n\to\infty}t_n=T_m,\qquad \sup_n u(t_n,0)<\infty.
\end{equation}
Let $t_*\in(0,T_m)$.
It follows from Lemma~\ref{Lemma:2.1}, \eqref{eq:1.6}, and \eqref{eq:2.35} that 
\begin{equation}
\label{eq:2.36}
\int_0^{s(t_n)}\,|\partial_x u(t_n)|^2\,\dee x +\frac{1}{2}\int_{t_*}^{t_n}\,(s')^3\,\dee\tau
 \le E(s(t_*),u(t_*))+\frac{1}{p+1}u(t_n,0)^{p+1}\le C
\end{equation}
for sufficiently large $n$. 
Then, by H\"older's inequality and \eqref{eq:2.36}, we see that
\begin{equation}
\label{eq:2.37}
s(t_n)=\int_{t_*}^{t_n}\, s'\,\dee\tau+s(t_*)
\le t_n^{\frac{2}{3}}\left(\int_{t_*}^{t_n}\, (s')^3\,\dee\tau\right)^{\frac{1}{3}}+s(t_*)\le C
\end{equation}
for sufficiently large $n$. 
By \eqref{eq:2.36} and \eqref{eq:2.37}, we obtain
$$
u(t_n,x)=u(t_n,x)-u(t_n,s(t_n))=-\int_x^{s(t_n)}\partial_x u(t_n)\,\dee x
\le s(t_n)^{\frac{1}{2}}\left(\int_0^{s(t_n)}|\partial_x u(t_n)|^2\,\dee x\right)^{\frac{1}{2}}\le C
$$
for $x\in [0,s(t_n)]$ and sufficiently large $n$.  
This, together with \eqref{eq:2.36} and \eqref{eq:2.37}, implies that 
$$
\sup_n\|u(t_n)\|_{W^{1,2}(0,s(t_n))}<\infty.
$$
Since $t_n \to T_m$ as $n \to \infty$ and $(s(t_n),u(t_n))\in W$, 
Propositions~\ref{Proposition:2.1} and \ref{Proposition:2.3} imply that 
the solution $(s,u)$ can be uniquely extended to a solution to problem~\eqref{SP} on $(0,T]$ for some $T > T_m$.
This contradicts the definition of $T_m$.
Hence, \eqref{eq:2.34} holds, and the proof of Proposition~\ref{Proposition:2.4} is complete.
$\Box$
\vspace{5pt}

Next, we establish the comparison principle for problem~\eqref{SP}.
\begin{proposition}
\label{Proposition:2.5}
Let $(s,u)$ and $(\sigma,v)$ be solutions to problem~\eqref{SP} with initial data $(s_0,\varphi)\in W$ and $(\sigma_0,\psi)\in W$, respectively.
Denote by $T_u$ and $T_v$ the maximal existence times of the solutions $(s,u)$ and $(\sigma,v)$, respectively.
Assume that 
$$
s_0\le\sigma_0,\qquad \varphi\le\psi\quad\mbox{on}\quad [0,s_0]. 
$$
Then 
\begin{equation}
\label{eq:2.38}
T_u\ge T_v,\qquad
\mbox{$s(t)\le \sigma(t)$ on $(0,T_v)$},\qquad 
\mbox{$u\le v$ in $D_s(T)$ for any $T\in(0,T_v)$}. 
\end{equation}
\end{proposition}
{\bf Proof.}
Consider the case where
\begin{equation}
\label{eq:2.39}
(s_0,\varphi)\in W,\quad (\sigma_0,\psi)\in W,\quad 
s_0<\sigma_0,\quad 
\varphi<\psi\mbox{ on $[0,s_0]$}. 
\end{equation}
Let $0<T<\min\{T_u,T_v\}$, and define
$$
T_*:=\sup\left\{t\in(0,T]\,:\,s(\tau)<\sigma(\tau)\mbox{ for $\tau\in(0,t)$}\right\}.
$$
It follows from \eqref{eq:2.39} that $T_*>0$ and 
\begin{equation}
\label{eq:2.40}
s(t)<\sigma(t)\quad\mbox{for}\quad t\in[0,T_*),\qquad s(T_*)=\sigma(T_*)\quad \mbox{if} \quad T_*<T.
\end{equation} 
Set 
$$
w(t,x):=v(t,x)-u(t,x),\quad (t,x)\in D_s(T_*). 
$$
Then, by \eqref{eq:2.39} and \eqref{eq:2.40}, $w$ satisfies
\begin{equation*}
\left\{
\begin{array}{ll}
\partial_t w=\partial_x^2 w, & (t,x)\in D_s(T_*),\vspace{3pt}\\
-\partial_xw(t,0)=a(t)w(t,0), & t\in(0,T_*],\vspace{3pt}\\
w(t,s(t))=v(t,s(t))>0, & t\in(0,T_*),\vspace{3pt}\\
w(T_*,s(T_*))=0 & \mbox{if $T_*<T$},\vspace{3pt}\\
w(0,x)>0, & x\in[0,s_0],
\end{array}
\right.
\end{equation*}
where $a$ is defined as in \eqref{eq:2.5}. 
Then, similarly to \eqref{eq:2.6}, we obtain 
\begin{equation}
\label{eq:2.41}
w>0\quad\mbox{in}\quad D_s^{**}(T_*)\setminus\{(T_*,s(T_*))\}. 
\end{equation}
Indeed, if $w$ is nonpositive for some point $(t_0,x_0)\in D_s^{**}(T_*)\setminus\{(T_*,s(T_*))\}$, then, 
by the continuity of $w$, we can find $(t_1,x_1)\in D_s^*(T_*)\setminus\{(T_*,s(T_*))\}$ such that  
$$
w(t_1,x_1)=0,\qquad \mbox{$w(t,x)>0$ for $(t,x)\in D_s^{**}(T_*)$ with $t<t_1$}.
$$
The maximum principle implies that $x_1=0$ and $w>0$ in $D_s(t_1)$, 
which, together with Hopf's lemma, implies that $0>-\partial_x w(t_1,0)=a(t)w(t_1,0)=0$. 
This is a contradiction, and hence \eqref{eq:2.41} holds. 

Assume that $T_*<T$. 
By \eqref{eq:2.41}, we apply Hopf's lemma again to obtain
\begin{equation}
\label{eq:2.42}
0>\partial_x w(T_*,s(T_*))=\partial_x v(T_*,s(T_*))-\partial_x u(T_*,\sigma(T_*))
=-\sigma'(T_*)+s'(T_*).
\end{equation}
On the other hand, it follows from \eqref{eq:2.40} that $\sigma'(T_*)\le s'(T_*)$,
which contradicts \eqref{eq:2.42}.
Thus, $T_*=T$. 
Combining \eqref{eq:2.40} and \eqref{eq:2.41}, we then obtain
$$
s(t)<\sigma(t)\quad\mbox{for}\quad t\in[0,T),\qquad
u\le v\quad\mbox{in}\quad D_s(T),
$$
for any $T\in(0,\min\{T_u,T_v\})$. 
This, together with Proposition~\ref{Proposition:2.4}, implies that $T_u\ge T_v$.
Hence, \eqref{eq:2.38} holds, and Proposition~\ref{Proposition:2.5} follows in the case of \eqref{eq:2.39}.

We complete the proof of Proposition~\ref{Proposition:2.5}. 
Assume that $s_0\le\sigma_0$, $(s_0,\varphi)\in W$, $(\sigma_0,\psi)\in W$, and $\varphi\le\psi$ on $[0,s_0]$. 
We identify $\psi$ with its zero extension. 
For any $n=1,2,\dots$, let $(\sigma_n,v_n)$ be a solution to problem~\eqref{SP} with initial data $(\sigma_{0,n},\psi_n)\in W$, 
where
$$
\sigma_{0,n}:=\sigma_0+\frac{1}{n},\qquad 
\psi_n(x):=\psi(x)+\frac{1}{n}(\sigma_{0,n}-x)\quad\mbox{for}\quad x\in[0,\sigma_{0,n}].
$$
Denote by $T_{v_n}$ its maximal existence time of the solution~$(\sigma_n,v_n)$. 
By Proposition~\ref{Proposition:2.5} with \eqref{eq:2.39}, we have 
$$
T_u\ge T_{v_n},\qquad
\mbox{$s(t)\le \sigma_n(t)$ on $(0,T_{v_n})$},\qquad 
\mbox{$u\le v_n$ in $D_s(T)$ for any $T\in(0,T_{v_n})$}. 
$$
Passing to the limit as $n\to\infty$, and applying Propositions~\ref{Proposition:2.2} and \ref{Proposition:2.4}, 
we obtain 
$$
T_u\ge T_v,\qquad
\mbox{$s(t)\le \sigma(t)$ on $(0,T_v)$},\qquad 
\mbox{$u\le v$ in $D_s(T)$ for any $T\in(0,T_v)$}. 
$$
Thus, Proposition~\ref{Proposition:2.5} follows.
$\Box$\vspace{7pt}

By Propositions~\ref{Proposition:2.4} and \ref{Proposition:2.5}, 
we improve Proposition~\ref{Proposition:2.2} to obtain the following result. 
\begin{proposition}
\label{Proposition:2.6}
Let $(s,u)$ be a solution to problem~\eqref{SP} with initial data $(s_0,\varphi)\in W$.
Let $T_m\in(0,\infty]$ denote the maximal existence time of $(s,u)$. 
Then, for any $T\in(0,T_m)$ and any $\epsilon>0$, 
there exists $\eta>0$ with the following property:
\begin{itemize}
  \item If $(\sigma_0,\psi)\in W$ satisfies
  $$
  |s_0-\sigma_0|+\|\varphi-\psi\|_{L^\infty}<\eta, 
  $$
  then problem~\eqref{SP} admits a solution~$(\sigma,v)$ on $(0,T]$ with initial data $(\sigma_0,\psi)$ 
  such that 
  $$
  \sup_{t\in(0,T]}|s(t)-\sigma(t)|+\sup_{t\in(0,T]}\|u(t)-v(t)\|_{L^\infty}<\epsilon.
  $$
\end{itemize}
\end{proposition}
{\bf Proof.}
By Propositions~\ref{Proposition:2.2}, \ref{Proposition:2.4}, and \ref{Proposition:2.5}, 
we may apply the same argument as in the proof of \cite{S}*{Theorem~2.2} 
to conclude Proposition~\ref{Proposition:2.6}. 
We therefore omit the details.
$\Box$
%%%%%%%%%%%%%%%%%%%%%%%%%%%%%%
%%%%%%%%%%%%%%%%%%%%%%%%%%%%%%
\section{Finite-time blow-up solutions}
%%%%%%%%%%%%%%%%%%%%%%%%%%%%%%
%%%%%%%%%%%%%%%%%%%%%%%%%%%%%%
We first give a sufficient condition for solutions to problem~\eqref{SP} to blow up in finite time. 
\begin{proposition}
\label{Proposition:3.1}
Let $(s,u)$ be a solution to problem~\eqref{SP} with initial data $(s_0,\varphi)\in W$, 
and denote by $T_m$ its maximal existence time. 
Assume that 
\begin{equation}
\label{eq:3.1}
u\in C^{0,1}(D_s^{**}(T))\quad\mbox{for some $T\in(0,T_m)$}. 
\end{equation}
If 
\begin{equation}
\label{eq:3.2}
E(s_0,\varphi):=\frac{1}{2}\int_0^{s_0} |\partial_x\varphi|^2\,\dee x-\frac{1}{p+1}\varphi(0)^{p+1}
<\frac{\pi^2}{256}\frac{\|\varphi\|_{L^1}^3}{(s_0+\|\varphi\|_{L^1})^4},
\end{equation}
then the solution~$(s,u)$ blows up in finite time, that is, $T_m<\infty$.
\end{proposition}
For the proof, we prepare the following lemma.
\begin{lemma}
\label{Lemma:3.1}
Let $(s,u)$ be a global-in-time solution to problem~\eqref{SP} with initial data $(s_0,\varphi)\in W$.
Then
$$
\int_0^\infty (s')^3\,\dee t\ge\frac{1}{128}\frac{\pi^2\|\varphi\|_{L^1}^3}{(s_0+\|\varphi\|_{L^1})^4}.
$$
\end{lemma}
{\bf Proof.}
Let $(s,u)$ be a global-in-time solution to problem~\eqref{SP} with initial data $(s_0,\varphi)\in W$.
Consider the free boundary problem
\begin{equation*}
\left\{
\begin{array}{ll}
\partial_t v=\partial_x^2 v, & (t,x)\in D_\sigma,\vspace{3pt}\\
\partial_x v(t,0)=0, & t\in(0,\infty),\vspace{3pt}\\
v(t,\sigma(t))=0, & t\in(0,\infty),\vspace{3pt}\\
\sigma'(t)=-v_x(t,\sigma(t)), & t\in(0,\infty),\vspace{3pt}\\
v(0,x)=\varphi(x), & x\in[0,s_0],\vspace{3pt}\\
\sigma(0)=s_0.
\end{array}
\right.
\end{equation*}
By \cite{Fr}*{Theorem~2, Chapter 8}, there exists a solution $(\sigma,v)$ to this problem.
We apply arguments similar to those in the proof of Proposition~\ref{Proposition:2.5} to obtain
\begin{align}
\label{eq:3.3}
  & s(t)\ge\sigma(t)\ge s_0,\quad t\in(0,\infty),\\
\nonumber
 & u(t,x)\ge v(t,x),\quad (t,x)\in D_\sigma.
\end{align}
Furthermore, it follows that  
\begin{align*}
\frac{\dee}{\dee t}\int_0^{\sigma(t)}v(t)\,\dee x
 & =\int_0^{\sigma(t)}\partial_tv(t)\,\dee x+v(t,\sigma(t))\sigma'(t)=\int_0^{\sigma(t)}\partial_x^2 v(t)\,\dee x\\
 & =\partial_xv(t,\sigma(t))-\partial_xv(t,0)=-\sigma'(t),\quad t\in(0,\infty). 
\end{align*}
This, together with $\sigma(0)=s_0$, implies that  
\begin{equation}
\label{eq:3.4}
\sigma(t)-s_0=\|\varphi\|_{L^1}-\|v(t)\|_{L^1},\quad t\in(0,\infty).
\end{equation} 
On the other hand, we observe from H\"older's inequality that 
$$
\int_0^\infty (s')^3\,\dee t\ge \int_0^t (s')^3\,\dee t
\ge t^{-2}\left(\int_0^t s'\, \dee t \right)^3=t^{-2}(s(t)-s_0)^3,\quad t\in(0,\infty),
$$
which, together with \eqref{eq:3.3} and \eqref{eq:3.4}, 
implies that
\begin{equation}
\label{eq:3.5}
\begin{split}
\int_0^\infty (s')^3\,\dee t \ge t^{-2}(\sigma(t)-s_0)^3
= t^{-2}(\|\varphi\|_{L^1}-\|v(t)\|_{L^1})^3,\quad t\in(0,\infty). 
\end{split}
\end{equation}

Define
$$
w(t,x):=(4\pi t)^{-\frac{1}{2}}\int_0^\infty
\left\{\exp\left(-\frac{|x-y|^2}{4t}\right)+\exp\left(-\frac{|x+y|^2}{4t}\right)\right\}\varphi(y)\,\dee y
$$
for $(t,x)\in(0,\infty)\times(0,\infty)$. 
Since $w$ is a solution to the heat equation in $(0,\infty)\times(0,\infty)$ with the zero Neumann boundary condition, 
the comparison principle yields
\[
\|v(t)\|_{L^\infty}\le\|w(t)\|_{L^\infty}\le(\pi t )^{-\frac{1}{2}}\|\varphi\|_{L^1},\quad t\in(0,\infty),
\]
which, together with \eqref{eq:3.4}, implies that 
$$
\|v(t)\|_{L^1}\le \sigma(t)\|v(t)\|_{L^\infty}
\le\sigma(t)(\pi t)^{-\frac{1}{2}}\|\varphi\|_{L^1}
\le(s_0+\|\varphi\|_{L^1})(\pi t)^{-\frac{1}{2}}\|\varphi\|_{L^1}
$$
for $t\in(0,\infty)$.
Hence, 
$$
\|v(t_0)\|_{L^1}\le\frac{1}{2}\|\varphi\|_{L^1},
\quad\mbox{where}\quad t_0:=4\pi^{-1}(s_0+\|\varphi\|_{L^1})^2.
$$
Then, applying \eqref{eq:3.5} with $t=t_0$, we obtain 
$$
\int_0^\infty (s')^3\,\dee t
\ge\left(4\pi^{-1}(s_0+\|\varphi\|_{L^1})^2\right)^{-2} \left(\|\varphi\|_{L^1}-\frac{1}{2}\|\varphi\|_{L^1}\right)^3
=\frac{\pi^2}{128}\frac{\|\varphi\|_{L^1}^3}{(s_0+\|\varphi\|_{L^1})^4}.
$$
This completes the proof of Lemma~\ref{Lemma:3.1}. 
$\Box$\vspace{5pt}

We are ready to prove Proposition~\ref{Proposition:3.1}.
\vspace{5pt}
\newline
{\bf Proof of Proposition~\ref{Proposition:3.1}.} 
Let $(s,u)$ be a solution to problem~\eqref{SP} with initial data $(s_0,\varphi)\in W$, and satisfy \eqref{eq:3.1} and \eqref{eq:3.2}. 
Assume that the solution~$(s,u)$ exists globally in time. 
By \eqref{eq:3.1}, we have
\begin{equation}
\label{eq:3.6}
\lim_{t\to+0}E(s(t),u(t))=E(s_0,\varphi).
\end{equation}
Define a function $F$ on $[0,\infty)$ by
\begin{equation}
\label{eq:3.7}
F(t)=\int_0^t\int_0^{s(\tau)}u^2\,\dee x\,\dee\tau,\quad t\in[0,\infty).
\end{equation}
It follows that
\begin{equation}
\label{eq:3.8}
\begin{split}
F''(t) &=\frac{\dee}{\dee t}\int_0^{s(t)}u(t)^2\,\dee x=2\int_0^{s(t)}u(t)\partial_tu(t)\,\dee x+u(t,s(t))^2s'(t)\\
 & =2\int_0^{s(t)}u(t)\partial_tu(t)\,\dee x=2\int_0^{s(t)}u(t)\partial_x^2u(t)\,\dee x\\
 & =-2\int_0^{s(t)}|\partial_xu(t)|^2\,\dee x+2u(t,s(t))\partial_xu(t,s(t))-2u(t,0)\partial_xu(t,0)\\
& =-2\int_0^{s(t)}|\partial_xu(t)|^2\,\dee x+2u(t,0)^{p+1}\\
&=-2(p+1)E(s(t),u(t))+(p-1)\int_0^{s(t)}|\partial_xu(t)|^2\,\dee x
\end{split}
\end{equation}
for $t\in(0,\infty)$.
This, together with \eqref{eq:2.8} and \eqref{eq:3.6},  implies that
\begin{equation}
\label{eq:3.9}
\begin{split}
F''(t) & =-2(p+1)\left(\int_0^t \frac{\dee}{\dee t}E(s(\tau),u(\tau))\,\dee\tau
+E(s_0,\varphi) \right)+(p-1)\int_0^{s(t)}|\partial_xu(t)|^2\,\dee x\\
 & =2(p+1)\int_0^t\int_0^{s(\tau)}|\partial_tu|^2\,\dee x\,\dee\tau\\
 & \qquad
 +2(p+1)\left(\frac{1}{2}\int_{t_0}^t (s')^3\,\dee\tau-E(s_0,\varphi)\right)+(p-1)\int_0^{s(t)}|\partial_xu(t)|^2\,\dee x
\end{split}
\end{equation}
for $t\in(0,\infty)$. 
Since $(s,u)$ is a global-in-time solution, 
it follows from Lemma~\ref{Lemma:3.1} and \eqref{eq:3.2} that
$$
E(s_0,\varphi)
<\frac{1}{2}\frac{\pi^2}{128}\frac{\|\varphi\|_{L^1}^3}{(s_0+\|\varphi\|_{L^1})^4}
\le\frac{1}{2}\int_0^\infty (s')^3\,\dee\tau.
$$ 
This implies that 
$$
E(s_0,\varphi)
\le\frac{1}{2}\int_0^t (s')^3\,\dee\tau
$$
for sufficiently large $t$.
Then, by \eqref{eq:3.9}, we have
\begin{equation}
\label{eq:3.10}
F''(t)>2(p+1)\int_0^t\int_0^{s(\tau)}|\partial_tu|^2\,\dee x\,\dee\tau>0
\end{equation}
for sufficiently large $t$.
Therefore, 
by H\"older's inequality, \eqref{eq:3.7}, and \eqref{eq:3.8}, we obtain
\begin{equation}
\label{eq:3.11}
\begin{split}
FF''(t) &\ge2(p+1)\int_0^t\int_0^{s(\tau)}|\partial_t u|^2\,\dee x\,\dee\tau
\int_0^t\int_0^{s(\tau)}u^2\,\dee x\,\dee\tau\\
 &\ge 2(p+1)\left(\int_0^t\int_0^{s(\tau)}u \partial_t u\,\dee x\,\dee\tau\right)^2
 =2(p+1)\left(\frac{1}{2}\int_0^t F''(\tau)\,\dee\tau\right)^2\\
 &  =\frac{p+1}{2}(F'(t)-F'(0))^2
\end{split}
\end{equation}
for sufficiently large $t$.

On the other hand, by \eqref{eq:3.10}, 
there exists $t_1\in(0,\infty)$ such that 
$$
F'(t)\ge F'(t_1)=\int_0^{s(t_1)}u(t_1)^2\,\dee x>0,\quad t\in[t_1,\infty),
$$
which implies that $\lim_{t\to\infty}F(t)=\infty$. 
It follows from \eqref{eq:3.10} that
$$
\lim_{t\to\infty}F'(t)=\infty.
$$
Then we have 
\begin{equation}
\label{eq:3.12}
\lim_{t\to\infty}\frac{(F'(t)-F'(0))^2}{(F'(t))^2}=1.
\end{equation}
Consequently, by \eqref{eq:3.11} and \eqref{eq:3.12}, we obtain
\begin{equation}
\label{eq:3.13}
\frac{F(t)F''(t)}{F'(t)^2}\ge\frac{p+1}{2}\frac{(F'(t)-F'(0))^2}{(F'(t))^2}
>\frac{p+3}{4}
\end{equation}
for sufficiently large $t$. 
Here, we used the inequality 
$p+3<2(p+1)$. 

Setting 
$$
G(t):=F(t)^{-\alpha}\quad\mbox{with}\quad \alpha=\frac{p-1}{4},
$$
we see that 
\begin{equation}
\label{eq:3.14}
G(t)>0,\quad
G'(t)=-\alpha F'(t)F(t)^{-\alpha-1}<0,
\end{equation}
for sufficiently large $t$. 
Furthermore, by \eqref{eq:3.13}, we obtain
\begin{equation}
\label{eq:3.15}
G''(t)=-\alpha F(t)^{-\alpha-2}
\left(F(t)F''(t)-\frac{p+3}{4}F'(t)^2\right)<0
\end{equation}
for sufficiently large $t$. 
By \eqref{eq:3.14} and \eqref{eq:3.15},
we conclude that the function~$G$ is strictly concave, strictly decreasing, and positive for all sufficiently large $t$. 
This leads to a contradiction. 
Thus, Proposition~\ref{Proposition:3.1} follows.
$\Box$\vspace{5pt}

As a corollary of Proposition~\ref{Proposition:3.1}, we have:
\begin{proposition}
\label{Proposition:3.2}
Let $(s,u)$ be a global-in-time solution to problem~\eqref{SP} with initial data $(s_0,\varphi)\in W$. 
Then 
$$
\inf_{t>0}E(s(t),u(t))\ge 0,\quad \int_0^\infty (s')^3\,\dee t<\infty.
$$
\end{proposition}
{\bf Proof.} 
Proposition~\ref{Proposition:3.1} implies that
if $E(s(t),u(t))<0$ for some $t>0$, then the solution $(s,u)$ blows up in finite time.
Hence, 
$$
\inf_{t>0}E(s(t),u(t))\ge 0.
$$ 
Then, by \eqref{eq:2.10} and \eqref{eq:2.13}, we obtain 
$$
\int_0^\infty (s')^3\,\dee t\le 2E(s(1),u(1))+\int_0^1 (s')^3\,\dee t<\infty.
$$
Thus, Proposition~\ref{Proposition:3.2} follows.
$\Box$\vspace{5pt}

In the rest of this section, we describe the behavior of blow-up solutions at the blow-up time.
The blow-up of the energy is discussed separately in Proposition~\ref{Proposition:4.1}.
\begin{proposition}
\label{Proposition:3.3}
Let $(s,u)$ be a solution to problem~\eqref{SP} with initial data $(s_0,\varphi) \in W$, 
and let $T_m$ denote its maximal existence time. 
If $T_m<\infty$, then \eqref{eq:1.7} and \eqref{eq:1.8} hold. 
Furthermore, 
\begin{equation}
\label{eq:3.16}
\lim_{t\nearrow T_m}s(t)<\infty.
\end{equation}
\end{proposition}
{\bf Proof.}
The proof of \eqref{eq:1.7} follows the arguments in the proof of \cite{FS2}*{Theorem~2.1},
with the aid of Liouville-type theorems. 
For any $t\in(0,T_m)$, define 
\begin{equation}
\label{eq:3.17}
M(t):=\sup_{\tau\in(0,t)}u(\tau,0),\qquad \lambda(t):=M(t)^{-(p-1)}.
\end{equation}
Then $M(t)$ is a positive, continuous, and nondecreasing function on $(0,T_m)$. 
Furthermore, it follows from Proposition~\ref{Proposition:2.4} that $M(t)\to\infty$ as $t\to T_m$. 
Consequently, for any $t\in(0,T_m)$, there exists $\xi(t)\in(0,T_m)$ such that 
\begin{equation}
\label{eq:3.18}
\xi(t):=\sup\{\xi\in(0,T_m):M(\xi)=2M(t)\}.
\end{equation}
It then suffices to show that there exists $C_*>0$ such that 
\begin{equation}
\label{eq:3.19}
\lambda(t)^{-2}(\xi(t)-t)\le C_*,\qquad t\in(T_m/2,T_m).
\end{equation}
Indeed, if \eqref{eq:3.19} holds, 
then for any $t\in(T_m/2,T_m)$, 
setting $t_0=t$ and $t_{n+1}=\xi(t_n)$ for $n=0,1,2,\dots$, 
we obtain  
\begin{align*}
t_n & \le t_{n-1}+C_*\lambda(t_{n-1})^2=t_{n-1}+C_*M(t_{n-1})^{-2(p-1)}\\
 & \le t_{n-2}+C_*M(t_{n-2})^{-2(p-1)}+C_*2^{-2(p-1)}M(t_{n-2})^{-2(p-1)}\\
 & \le t_0+C_*M(t_0)^{-2(p-1)}\sum_{k=0}^{n-1} 2^{-2k(p-1)}
\end{align*}
for $n=1,2,\dots$. Since $t_n\to T_m$ as $n\to\infty$, it follows that
$$
T_m\le t+CC_*M(t)^{-2(p-1)},
$$
which implies \eqref{eq:1.7}. 

We prove \eqref{eq:3.19} by contradiction. 
Assume that there exists a sequence $\{t_n\}\subset(T_m/2,T_m)$ such that 
\begin{equation}
\label{eq:3.20}
\lim_{j\to\infty}\lambda(t_n)^{-2}(\xi(t_n)-t_n)=\infty. 
\end{equation}
For any $n=1,2,\dots$, we take a sequence $\{\hat{t}_n\}\subset(0,t_n]$ satisfying 
\begin{equation}
\label{eq:3.21}
u(\hat{t}_n,0)\ge\frac{1}{2}M(t_n). 
\end{equation}
It follows from $\lim_{n\to\infty}t_n=T_m$ that $\lim_{n\to\infty}M(t_n)=\infty$, which implies that 
\begin{equation}
\label{eq:3.22}
\lim_{n\to\infty}\hat{t}_n=T_m.
\end{equation}
Set $\lambda_n:=\lambda(t_n)$, and define
$$
v_n(\xi,y):=\lambda_n^{\frac{1}{p-1}}u(\lambda_n^2\xi+\hat{t}_n,\lambda_n y)
$$
for $\xi\in(-\lambda_n^{-2}\hat{t}_n,\lambda_n^{-2}(T_m-\hat{t}_n))$ and $y\in(0,\lambda_n^{-1}s(\lambda_n^2\xi+\hat{t}_n))$. 
Then $v_n$ satisfies 
\begin{equation*}
\left\{
\begin{array}{ll}
\partial_\xi v_n=\partial_y^2 v_n,\quad & 
\xi\in(-\lambda_n^{-2}\hat{t}_n,\lambda_n^{-2}(T_m-\hat{t}_n)),\,\,\, y\in(0,\lambda_n^{-1}s(\lambda_n^2\xi+\hat{t}_n)),\vspace{5pt}\\
-\partial_yv(\xi,0)=v(\xi,0)^p, & \xi\in(-\lambda_n^{-2}\hat{t}_n,\lambda_n^{-2}(T_m-\hat{t}_n)),\vspace{5pt}\\
v_n(\xi,\lambda_n^{-1}s(\lambda_n^2\xi+\hat{t}_n))=0, & \xi\in(-\lambda_n^{-2}\hat{t}_n,\lambda_n^{-2}(T_m-\hat{t}_n)).
\end{array}
\right.
\end{equation*}
By the maximum principle, together with \eqref{eq:3.17} and \eqref{eq:3.18}, we obtain
\begin{equation}
\label{eq:3.23}
0\le v_n(\xi,y)\le \frac{1}{M(t_n)}\max\left\{\|\varphi\|_{L^\infty},\sup_{0\le t\le\xi(t_n)} u(t,0)\right\}\le 2
\end{equation}
for $\xi\in(-\lambda_n^{-2}\hat{t}_n,\lambda_n^{-2}(\xi(t_n)-\hat{t}_n))$, $y\in[0,\lambda_n^{-1}s(\lambda_n^2\xi+\hat{t}_n))$, and sufficiently large $n$. 
Furthermore, by \eqref{eq:1.3}, \eqref{eq:3.21}, and \eqref{eq:3.22}, we deduce that
\begin{equation}
\label{eq:3.24}
v_n(0,0)\ge\frac{1}{2}
\end{equation}
and 
\begin{equation}
\label{eq:3.25}
\lambda_n^{-1}s(\lambda_n^2\xi+\hat{t}_n)\ge \lambda_n^{-1}s_0\to\infty,\qquad
-\lambda_n^{-2}\hat{t}_n\to-\infty,
\end{equation}
as $n\to\infty$. 
Thanks to \eqref{eq:3.20}, \eqref{eq:3.23}, and \eqref{eq:3.25},
by applying parabolic regularity theorems, the Arzel\'a--Ascoli theorem, and a diagonal argument 
and by passing to a subsequence if necessary, 
we obtain a function 
$$
v\in C^{1;2}_{\rm loc}((-\infty,\infty)\times(0,\infty))\cap C^{0;1}_{\rm loc}((-\infty,\infty)\times[0,\infty))
$$
such that 
\begin{align*}
 & \lim_{n\to\infty}\|v_n-v\|_{C^{1;2}(K)}=0\quad\mbox{for any compact set $K\subset(-\infty,\infty)\times(0,\infty)$}, \\
 & \lim_{n\to\infty}\|v_n-v\|_{C^{0;1}(K')}=0\quad\mbox{for any compact set $K'\subset(-\infty,\infty)\times[0,\infty)$}.
\end{align*}
These, together with \eqref{eq:3.23}, imply that 
the limit function $v$ satisfies
\begin{equation}
\label{eq:3.26}
\left\{
\begin{array}{ll}
\partial_\xi v=\partial_y^2 v, & (\xi,y)\in(-\infty,\infty)\times(0,\infty),\vspace{5pt}\\
 0\le v\le 2, & (\xi,y)\in(-\infty,\infty)\times(0,\infty),\vspace{5pt}\\
 -\partial_y v(\xi,0)=v(\xi,0)^p, & \xi\in(-\infty,\infty).
\end{array}
\right.
\end{equation}
In addition, it follows from \eqref{eq:3.24} that $v(0,0)\ge 1/2$.
However, applying \cite{Q1}*{Theorem~1} to \eqref{eq:3.26}, we conclude that $v$ must be identically zero in $(-\infty,\infty)\times[0,\infty)$, 
which contradicts the fact that $v(0,0)\ge 1/2$. 
Thus, \eqref{eq:3.19} holds, and hence we obtain \eqref{eq:1.7}.

We prove \eqref{eq:1.8} and \eqref{eq:3.16}. 
Let $\delta\in (0,s_0)$. 
Since $s(t)\ge s_0$ for $t\in(0,T_m)$, 
by \eqref{eq:1.7}, we apply \cite{HY}*{Theorem~4.1} to obtain 
$$
\sup_{t\in(0,T_m)}\sup_{x\in[\delta/2,s_0/2]}u(t,x)<\infty.
$$
Then we apply parabolic regularity theorems to obtain 
$$
|\partial_x u(t,\delta)|+|\partial_t u(t,\delta)|\le C
$$
for $t\in(T_m/2,T_m)$. 
This, together with \eqref{eq:2.11}, implies that
\begin{equation*}
\begin{split}
 & \frac{\dee}{\dee t}\int_\delta^{s(t)}|\partial_x u(t)|^2\,\dee x
=2\int_\delta^{s(t)} \partial_x u(t)\partial_x\partial_t u(t)\,\dee x+|\partial_x u(t,s(t))|^2s'(t)\\
 & =2\partial_x u(t,s(t))\partial_t u(t,s(t))-2\partial_x u(t,\delta)\partial_t u(t,\delta)
 -2\int_\delta^{s(t)} |\partial_tu(t)|^2\,\dee x+s'(t)^3\\
 & \le -s'(t)^3+C
\end{split}
\end{equation*}
for $t\in(T_m/2,T_m)$. It follows that
\begin{equation}
\label{eq:3.27}
\sup_{t\in(T_m/2,T_m)}\int_\delta^{s(t)}|\partial_x u(t)|^2\,\dee x+\int_{T_m/2}^{T_m} (s')^3\,\dee t\le CT_m+\int_\delta^{s(T_m/2)}\left|\partial_x u\left(\frac{T_m}{2}\right)\right|^2\,\dee x<\infty.
\end{equation}
Hence,
\begin{equation}
\label{eq:3.28}
\lim_{t\nearrow T_m}s(t)=s\left(\frac{T_m}{2}\right)+\int_{T_m/2}^{T_m} s'\,\dee t
\le s\left(\frac{T_m}{2}\right)+\left(\frac{T_m}{2}\right)^{\frac{2}{3}}\left(\int_{T_m/2}^{T_m} (s')^3\,\dee t\right)^{\frac{1}{3}}<\infty, 
\end{equation}
which is exactly~\eqref{eq:3.16}.
Furthermore, 
by \eqref{eq:3.27} and \eqref{eq:3.28}, we obtain 
$$
u(t,x)=u(t,x)-u(t,s(t))=-\int_x^{s(t)}\partial_x u(t)\,\dee x
\le s(t)^{\frac{1}{2}}\left(\int_x^{s(t)}|\partial_x u(t)|^2\,\dee x\right)^{\frac{1}{2}}\le C
$$
for $t\in(T_m/2,T_m)$ and $x\in(\delta,s(t))$. This implies \eqref{eq:1.8}. 
Therefore, the proof of Proposition~\ref{Proposition:3.3} is complete.
$\Box$
%%%%%%%%%%%%%%%%%%%%%%%%%%%%%%
%%%%%%%%%%%%%%%%%%%%%%%%%%%%%%
\section{Decay properties of global-in-time solutions}
%%%%%%%%%%%%%%%%%%%%%%%%%%%%%%
%%%%%%%%%%%%%%%%%%%%%%%%%%%%%%
In this section, we establish several propositions on the decay properties of global-in-time solutions to problem~\eqref{SP}. 
We first show that for any global-in-time solution $(s,u)$, 
the $L^\infty$-norm $\|u(t)\|_{L^\infty}$ converges to zero as $t \to \infty$.
\begin{proposition}
\label{Proposition:4.1}
Let $(s,u)$ be a global-in-time solution to problem~\eqref{SP} with initial data $(s_0,\varphi)\in W$. 
Then 
$$
\lim_{t \to \infty} \|u(t)\|_{L^\infty}=0.
$$
\end{proposition}
The following lemma is a key ingredient in the proof of Proposition~\ref{Proposition:4.1}. 
\begin{lemma}
\label{Lemma:4.1}
Let $(s,u)$ be a global-in-time solution to problem~\eqref{SP} with initial data $(s_0,\varphi)\in W$. 
Then $\lim_{t \to \infty}u(t,0)=0$.
\end{lemma}
{\bf Proof.} 
The proof is by contradiction. 
Let $(s,u)$ be a global-in-time solution to problem~\eqref{SP} with initial data $(s_0,\varphi)\in W$. 
Assume that
$$
\ell:=\limsup_{t\to \infty}\, u(t,0)\in(0,\infty].
$$
Then there exists a sequence $\{t_n\}_{n=0}^\infty\subset(0,\infty)$ such that
\begin{equation}
\label{eq:4.1}
t_0<t_1<\cdots<t_n<\cdots,\quad \lim_{n\to\infty}t_n=\infty,
\end{equation}
and 
\begin{itemize}
\item[(A)]
    if $\ell<\infty$, then
   \begin{equation}
    \label{eq:4.2}
    \sup_{t>t_0} u(t,0)\le\frac{3}{2}\ell,\quad 
    \frac{\ell}{2}\le u(t_n,0)\le\frac{3}{2}\ell\quad\mbox{for $n=0,1,2,\dots$};
    \end{equation} 
\item[(B)] 
    if $\ell=\infty$, then
    \begin{equation}
    \label{eq:4.3}
    \lim_{n\to\infty}u(t_n,0)=\infty,\quad 
    \sup_{t\in[t_0,t_n]} u(t,0)=u(t_n,0)\ge 1\quad\mbox{for $n=1,2,\dots$}.
    \end{equation}
\end{itemize}
In the proof, generic constants $C$ are independent of $n$. 
For any $n=1,2,\dots$, set 
\begin{equation}
\label{eq:4.4}
\begin{split}
 & \sigma_n:=u(t_n,0),\quad
\lambda_n:=\sigma_n^{-(p-1)},\quad \tau_n:=\lambda_n^{-2}(t_n-t_0),\\
 & s_n(\tau):=\lambda_n^{-1}s(t_n+\lambda_n^2\tau)\quad\mbox{for}\quad\tau\in[-\tau_n,0].
\end{split}
\end{equation}
Then it follows from \eqref{eq:4.1}--\eqref{eq:4.4} that 
$\{\lambda_n\}$ is bounded and 
$\lim_{n\to\infty}\tau_n=\infty$.
Set
\begin{equation}
\label{eq:4.5} 
D_n:=\bigcup_{\tau\in(-\tau_n,0]}\{\tau\}\times[0,s_n(\tau)),
\quad
D'_n:=(-\tau_n,0]\times[0,\infty).
\end{equation}
Extending $u(t,x)$ by zero for $x \in (s(t),\infty)$,
we define $v_n\in H^1_{\rm loc}(D'_n)$ by
$$
v_n(\tau,y):=\lambda_n^{\frac{1}{p-1}}u(t_n+\lambda_n^2\tau,\lambda_ny),\quad (\tau,y)\in D'_n.
$$
It follows from the maximal principle that
$$
 v_n(\tau,y)=\frac{u(t_n+\lambda_n^2\tau,\lambda_ny)}{u(t_n,0)}\le \frac{\max\{\sup_{[t_0,t_n]}u(t,0),\|u(t_0)\|_{L^\infty(0,s(t_0))}\}}{u(t_n,0)}
$$
for $(\tau,y)\in D_n'$, 
which, together with \eqref{eq:4.2} and \eqref{eq:4.3}, yields
\begin{equation}
\label{eq:4.6}
v_n(0,0)=1, \quad 0\le v_n\le C \quad in \quad D'_n,
\end{equation}
for $n=1,2,\dots$.
Furthermore, $v_n$ satisfies
\begin{equation}
\left\{
\label{eq:4.7}
\begin{array}{ll}
 \partial_{\tau}v_n=\partial_{y}^2v_n, & \quad (\tau,y)\in D_n,\vspace{7pt}\\
 -\partial_y v_n(\tau,0)=v_n(\tau,0)^p, & \quad \tau\in(-\tau_n,0],\vspace{7pt}\\
v_n(\tau,s_n(\tau))=0, & \quad \tau\in(-\tau_n,0],\vspace{3pt}\\
-\partial_y v_n(\tau,s_n(\tau))=\lambda_n^{\frac{1}{p-1}}s_n'(\tau), & \quad \tau\in(-\tau_n,0],
\end{array}
\right.
\end{equation}
for $n=1,2\dots$. 
In particular, it follows from \eqref{eq:4.6} and \eqref{eq:4.7} that
\begin{equation}
\label{eq:4.8}
0\le -\partial_y v_n(\tau,0)\le C^p,\quad \tau\in[-\tau_n,0],
\end{equation}
for $n=1,2,\dots$. 

By the boundedness of  $\{\lambda_n\}$, there exists $R>0$ such that
$$
0<R<\liminf_{n\to\infty}s_n(-\tau_n)=s(t_0)\liminf_{n\to\infty}\lambda_n^{-1}.
$$
Then, by applying parabolic regularity theorems, the Arzel\'a--Ascoli theorem, and a diagonal argument 
and by passing a subsequence if necessary, 
we obtain a function 
$$
v\in C^{1;2}((-\infty,0]\times[0,R))
$$ 
satisfying
\begin{equation}
\label{eq:4.9}
\lim_{n\to\infty}\|v_n-v\|_{C^{1;2}((-m,0]\times[0,R))}=0\quad\mbox{for $m>1$},
\end{equation}
and 
\begin{equation}
\left\{
\label{eq:4.10}
\begin{array}{ll}
 \partial_{\tau}v=\partial_{y}^2v, & \quad (\tau,y)\in(-\infty,0)\times(0,R),\vspace{3pt}\\
 -\partial_y v(\tau,0)=v(\tau,0)^p, & \quad \tau\in(-\infty,0],\vspace{3pt}\\
 v(0,0)=1.
\end{array}
\right.
\end{equation}
\vspace{5pt}
\newline
\underline{Step 1.}
We claim that for any $m\ge 1$, 
\begin{equation}
\label{eq:4.11}
\sup_{n\ge 0}\sup_{\tau\in[-m,0]} \,\int_0^m |\partial_y v_n(\tau)|^2\,\dee y
+\sup_{n\ge 0}\,\int_{-m}^0 \int_0^m
|\partial_\tau v_n|^2\,\dee y\,\dee \tau\le C.
\end{equation}
Let $\phi\in C_0^{\infty}([0,\infty))$ be such that 
\begin{equation}
\label{eq:4.12}
0\le\phi\le 1\mbox{ in $[0,\infty)$},\quad 
\phi=1\mbox{ in $[0,m]$},\quad
\phi=0\mbox{ in $[2m,\infty)$},\quad
|\partial_y\phi|\le\frac{2}{m}\mbox{ in $[0,\infty)$}.
\end{equation}
Recalling that $v_n(\tau,s_n(\tau))=0$ for $\tau\in(-\tau_n,0)$, 
by \eqref{eq:4.7} and \eqref{eq:4.12}, we obtain
\begin{align*}
 & \frac{1}{2}\frac{\dee}{\dee\tau}\int_0^{2m}v_n(\tau)^2\phi^2\,\dee y=\frac{1}{2}\frac{\dee}{\dee\tau}\int_0^{s_n(\tau)}v_n(\tau)^2\phi^2\,\dee y\\ 
 & =\int_0^{s_n(\tau)}\,v_n(\tau) \partial_\tau v_n(\tau)\phi^2\,\dee y
 =\int_0^{s_n(\tau)}\,v_n \partial_y^2v_n(\tau)\phi^2\, dy\\
 &=-v(\tau,0) \partial_yv_n(\tau,0)-\int_0^{s_n(\tau)}\,|\partial_y v_n(\tau)|^2\phi^2\,\dee y
 -2\int_0^{s_n(\tau)}\,v_n(\tau) \partial_yv_n(\tau) \phi \partial_y\phi\,\dee y\\
&\le -v(\tau,0) \partial_yv_n(\tau,0)
-\frac{1}{2}\int_0^{s_n(\tau)}\,|\partial_y v_n(\tau)|^2\phi^2\,\dee y
+C\int_0^{s_n(\tau)}\,v_n(\tau)^2|\partial_y\phi|^2\,\dee y
\end{align*}
for $\tau\in[-m,0]$ and sufficiently large $n$, 
that is, 
$$
\int_0^{s_n(\tau)}\,|\partial_y v_n(\tau)|^2\phi^2\,\dee y
 \le -2v_n(\tau,0)\partial_y v_n(\tau,0)
+C\int_0^{s_n(\tau)}\,v_n(\tau)^2 |\partial_y \phi|^2 \,\dee y
-\frac{\dee}{\dee\tau}\int_0^{2m}v_n(\tau)^2\phi^2\,\dee y
$$
for $\tau\in[-m,0]$ and sufficiently large $n$. 
Therefore, for any $m\ge 1$, 
by \eqref{eq:4.6}, \eqref{eq:4.7}, \eqref{eq:4.8}, and \eqref{eq:4.12},
we obtain
\begin{equation}
\label{eq:4.13}
\begin{split}
\int_{-m}^0\int_0^m\,|\partial_y v_n|^2\,\dee y\, d\tau
 & \le -2\int_{-m}^0\, v_n(\tau,0)\partial_y v_n(\tau,0)\,\dee\tau
 +C\int_{-m}^0 \int_0^{2m}\,v_n^2 \,\dee y\,d\tau\\
& \qquad -\int_0^m\,v_n(0,y)^2\,\dee y+\int_0^{2m}\,v_n(-m,y)^2\,\dee y\le C
\end{split}
\end{equation}
for sufficiently large $n$. 

On the other hand, it follows from \eqref{eq:4.7} that
$$
0=\frac{\dee}{\dee\tau}v_n(\tau,s_n(\tau))=\partial_\tau v_n(\tau,s_n(\tau))+\partial_y v_n(\tau,s_n(\tau))s_n'(\tau)
=\partial_\tau v_n(\tau,s_n(\tau))-\lambda_n^{\frac{1}{p-1}}\left(s_n'(\tau)\right)^2
$$
for $\tau\in[-m-1,0]$. 
Then it follows from \eqref{eq:4.6}--\eqref{eq:4.8} and \eqref{eq:4.12} that
\begin{equation*}
\begin{split}
& \frac{1}{2}\frac{\dee}{\dee\tau}\int_0^{s_n(\tau)}\,|\partial_y v_n(\tau)|^2\phi^2\,\dee y\\
& =\frac{1}{2}\lambda_n^{\frac{2}{p-1}}\left(s_n'(\tau)\right)^3\phi(s_n(\tau))^2
+\int_0^{s_n(\tau)}\,\partial_y v_n(\tau)\,\partial_\tau\partial_yv_n(\tau)\phi^2 \,\dee y\\\
& =\frac{1}{2}\lambda_n^{\frac{2}{p-1}}\left(s_n'(\tau)\right)^3\phi(s_n(\tau))^2+\partial_y v_n(\tau,s_n(\tau))\,\partial_\tau v_n(\tau,s_n(\tau))\phi(s_n(\tau))^2\\
& \qquad
-\partial_yv_n(\tau,0) \partial_\tau v_n(\tau,0)
-\int_0^{s_n(\tau)}\,\partial_y^2 v_n(\tau) \partial_\tau v_n(\tau)\phi^2\,\dee y
-2\int_0^{s_n(\tau)}\,\partial_y v_n(\tau) \partial_\tau v_n(\tau)\phi \partial_y\phi\,\dee y\\
& \le -\frac{1}{2}\lambda_n^{\frac{2}{p-1}}\left(s_n'(\tau)\right)^3\phi(s_n(\tau))^2
+v_n(\tau,0)^p\partial_\tau v_n(\tau,0)\\
 & \qquad
 -\int_0^{s_n(\tau)}\,|\partial_\tau v_n(\tau)|^2\phi^2\,\dee y+2\int_0^{s_n(\tau)}\,|\partial_y v_n(\tau)| |\partial_\tau v_n(\tau)| \phi |\partial_y\phi|\,\dee y\\
& \le \frac{1}{p+1}\partial_\tau \left(v_n(\tau,0)^{p+1}\right)-\frac{1}{2}\int_0^{s_n(\tau)}\,|\partial_\tau v_n(\tau)|^2\phi^2\,\dee y
+C\int_0^{s_n(\tau)}\,|\partial_y v_n(\tau)|^2 |\partial_y\phi|^2\,\dee y
\end{split}
\end{equation*}
for $\tau\in[-m-1,0]$ and sufficiently large $n$.
This, together with \eqref{eq:4.12}, implies that 
\begin{equation}
\label{eq:4.14}
\begin{split}
 & \frac{1}{2}\int_0^m\,|\partial_y v_n(\tau_2)|^2\,\dee y
 +\frac{1}{2}\int_{\tau_1}^{\tau_2}\int_0^{m}\,|\partial_\tau v_n|^2\,\dee y\,\dee\tau\\
 & \le \frac{1}{2}\int_0^{2m}\,|\partial_y v_n(\tau_1)|^2\,\dee y
 +\frac{1}{p+1}v_n(\tau_2,0)^{p+1}
 +C\int_{\tau_1}^{\tau_2}\int_0^{2m}\,|\partial_y v_n|^2\,\dee y\,\dee\tau
\end{split}
\end{equation}
for $\tau_1$, $\tau_2\in[-m-1,0]$ with $\tau_1<\tau_2$ and sufficiently large $n$. 
For any $\tau_2\in[-m,0]$, 
integrating both sides of \eqref{eq:4.14} in $\tau_1$ on $[-m-1,-m]$, 
we obtain 
\begin{equation}
\label{eq:4.15}
\begin{split}
\int_0^m\,|\partial_y v_n(\tau_2)|^2\,\dee y
 & \le \int_{-m-1}^{-m}\int_0^{2m}\,|\partial_y v_n|^2\,\dee y\,\dee\tau\\
 & +\frac{2}{p+1}v_n(\tau_2,0)^{p+1}+C\int_{-m-1}^{\tau_2}\int_0^{2m}\,|\partial_y v_n|^2\,\dee y\,\dee\tau. 
\end{split}
\end{equation}
Then, by \eqref{eq:4.6}, \eqref{eq:4.13} with $m$ replaced by $2m$, and \eqref{eq:4.15}, we have 
\begin{equation}
\label{eq:4.16}
\sup_{\tau\in[-m,0]}\int_0^m\,|\partial_y v_n(\tau)|^2\,\dee y\le C
\end{equation}
for sufficiently large $n$.
Furthermore, 
by \eqref{eq:4.6}, \eqref{eq:4.14} with $\tau_1=-m$ and $\tau_2=0$, and \eqref{eq:4.16}, 
we obtain 
$$
\sup_{\tau\in[-m,0]}\int_0^m\,|\partial_y v_n(\tau)|^2\,\dee y+\int_{-m}^0\int_0^m\,|\partial_\tau v_n|^2\,\dee y\,\dee\tau\le C
$$
for sufficiently large $n$. Thus, \eqref{eq:4.11} holds. 
By \eqref{eq:4.6} and \eqref{eq:4.11}, we conclude that
$$
\sup_n \|v_n\|_{H^1((-m,0)\times(0,m))}<\infty\quad\mbox{for}\quad m>1.
$$
Then, applying the weak convergence theorem in $H^1$, the Sobolev compact embedding theorem, and a diagonal argument, 
and passing to a subsequence if necessary, we obtain a function
\begin{equation}
\label{eq:4.17}
w\in H^1_{{\rm loc}}((-\infty,0]\times[0,\infty))
\end{equation}
such that 
\begin{equation}
\label{eq:4.18}
\begin{split}
\lim_{n\to\infty}v_n=w\quad & \mbox{weakly in $H^1(K)$ for any compact set $K \subset (-\infty,0]\times[0,\infty)$},\\
\lim_{n\to\infty}v_n=w\quad & \mbox{strongly in $L^2(K)$ for any compact set $K \subset (-\infty,0]\times[0,\infty)$},\\
\lim_{n\to\infty}v_n=w\quad & \mbox{for almost everywhere in $(-\infty,0)\times(0,\infty)$}.
\end{split}
\end{equation}
This, together with \eqref{eq:4.6} and \eqref{eq:4.9}, implies that
\begin{align}
\label{eq:4.19}
 & 0\le w\le C\quad\mbox{for almost everywhere in $(-\infty,0)\times(0,\infty)$},\\
\label{eq:4.20}
 & v=w \quad\mbox{for almost everywhere in $(-\infty,0)\times(0,R)$}.
\end{align}
\underline{Step 2.}
We complete the proof. 
Let $\zeta\in C_0^\infty((-\infty,0)\times(0,\infty))$. Let $m>1$ be such that 
$\mbox{supp}\,\zeta\subset Q:=(-m,0)\times(0,m)$. 
It follows from \eqref{eq:4.7} that
\begin{align*}
\int_0^\infty \partial_\tau v_n(\tau)\zeta(\tau)\,\dee y
 & =\int_0^{s_n(\tau)}\partial_\tau v_n(\tau)\zeta(\tau)\,\dee y=\int_0^{s_n(\tau)}\partial_y^2v_n(\tau)\zeta(\tau)\,\dee y\\
 & =\partial_y v_n(\tau,s_n(\tau))\zeta(\tau,s_n(\tau))-\int_0^{s_n(\tau)}\partial_y v_n(\tau)\partial_y\zeta(\tau)\,\dee y\\
&=-\lambda_n^{\frac{1}{p-1}}s_n'(\tau)\zeta(\tau,s_n(\tau))+\int_0^\infty v_n(\tau)\partial_y^2\zeta(\tau)\,\dee y
\end{align*}
for $\tau\in(-m,0)$. 
Then, by \eqref{eq:4.4}, we have
\begin{equation}
\begin{split}
\label{eq:4.21}
 & \left| \int_{-m}^0\int_0^\infty v_n\{-\partial_\tau\zeta-\partial_y^2\zeta\}\,\dee y\,\dee \tau\right|
 =\left| \int_{-m}^0\int_0^\infty \{\partial_\tau v_n\zeta-v_n\partial_y^2\zeta\}\,\dee y\,\dee \tau\right|\\
& =\lambda_n^{\frac{1}{p-1}}\left|\int_{-m}^0 \lambda_n s'(t_n+\lambda_n^2\tau)\zeta(\tau,s_n(\tau))\,\dee\tau\right| \\
& 
\le \lambda_n^{\frac{1}{p-1}-1}\|\zeta\|_{L^\infty(Q)} \int_{t_n-m\lambda_n^2}^{t_n} s'\,\dee t
\le \lambda_n^{\frac{1}{p-1}-1}\|\zeta\|_{L^\infty(Q)} \left(\int_{t_n-m\lambda_n^2}^{t_n} (s')^3\,\dee t \right)^{\frac{1}{3}}(m\lambda_n^2)^{\frac{2}{3}} \\
& \le \lambda_n^{\frac{1}{p-1}+\frac{1}{3}}\|\zeta\|_{L^\infty(Q)} m^{\frac{2}{3}} \left(\int_{t_n-m\lambda_n^2}^{t_n} (s')^3\,\dee t \right)^{\frac{1}{3}}.
\end{split}
\end{equation}

On the other hand, 
since $(s,u)$ is a global-in-time solution, 
it follows from Lemma~\ref{Lemma:2.1} and Proposition~\ref{Proposition:3.2} that 
\begin{equation}
\label{eq:4.22}
\int_{t_0}^\infty\int_0^{s(\tau)}|\partial_t u|^2\,\dee x\,\dee\tau+\frac{1}{2}\int_0^\infty (s')^3\,\dee\tau
\le E(s(t_0),u(t_0))<\infty. 
\end{equation}
This, together with the boundedness of $\{\lambda_n\}$, \eqref{eq:4.18}, and \eqref{eq:4.21}, implies that 
$$
0=\lim_{n\to\infty}\int_{-m}^0\int_0^\infty v_n\{-\partial_\tau\zeta-\partial_y^2\zeta\}\,\dee y\,\dee \tau
=\int_{-\infty}^0\int_0^\infty w\{-\partial_\tau\zeta-\partial_y^2\zeta\}\,\dee y\,\dee \tau, 
$$
that is,  
\begin{equation}
\label{eq:4.23}
\partial_\tau w=\partial_y^2 w\quad\mbox{in}\quad {\mathcal D}'((-\infty,0)\times(0,\infty)).
\end{equation}
Then, by \eqref{eq:4.17} and \eqref{eq:4.19}, we apply parabolic regularity theorems to deduce that
$$
w\in C^{1;2}((-\infty,0)\times(0,\infty))
$$ 
and that $w$ satisfies \eqref{eq:4.23} in the classical sense. 
Furthermore, by the boundedness of $\{\lambda_n\}$ and \eqref{eq:4.22}, 
we have
$$
\int_{-m}^0\int_0^m |\partial_\tau v_n|^2 \,\dee y \,\dee\tau\le\lambda_n^{\frac{p+1}{p-1}}\int_{t_n-m\lambda_n^2}^{t_n}\int_0^{s(t)}|\partial_t u|^2 \,\dee x\,\dee t \to 0
$$
as $n\to\infty$. This, together with \eqref{eq:4.18}, implies that $\partial_\tau w=0$ in $L^2_{\rm loc}((-\infty,0)\times(0,\infty))$. 
Hence, $w$ is independent of $\tau$.
It then follows from \eqref{eq:4.19} and \eqref{eq:4.23} that $w$ is a nonnegative bounded harmonic function in $(0,\infty)$.
Therefore, $w$ must be a nonnegative constant function in $(0,\infty)$.
By \eqref{eq:4.20}, we conclude that $v$ is a nonnegative constant function in $(-\infty,0)\times[0,R)$,
which contradicts the conditions $-\partial_y v(\tau,0)=v(\tau,0)^p$ for $\tau\in(-\infty,0]$ (see \eqref{eq:4.10}) and $v(0,0)=1$.
Hence, $\ell=0$, and Lemma~\ref{Lemma:4.1} follows.
$\Box$
\vspace{5pt}
\newline
{\bf Proof of Proposition~\ref{Proposition:4.1}.}
Let $(s,u)$ be a global-in-time solution to problem~\eqref{SP} with initial data $(s_0,\varphi)\in W$. 
Let $G_D$ be the Dirichlet heat kernel in $(0,\infty)\times(0,\infty)$, that is, 
$$
G_D(t,x,y):=(4\pi t)^{-\frac{1}{2}}\left(\exp\left(-\frac{|x-y|^2}{4t}\right)-\exp\left(-\frac{|x+y|^2}{4t}\right)\right)
$$
for $(t,x,y)\in(0,\infty)\times[0,\infty)^2$.
Set 
$$
z(t,x):=\int_0^\infty G_D(t,x,y)\varphi(y)\,\dee y+\int_0^t \partial_yG_D(t-\tau,x,0)u(\tau,0)\,\dee\tau
$$
for $(t,x)\in(0,\infty)\times[0,\infty)$. 
Then $z$ satisfies 
$$
\partial_t z=\partial_x^2 z\mbox{ in }(0,\infty)\times(0,\infty),
\quad
z(t,0)=u(t,0)\mbox{ for } t\in(0,\infty),
\quad
z(0,x)=\varphi(x)\mbox{ for } x\in(0,\infty).
$$
The comparison principle implies that 
$u\le z$ in $D_s$. 

On the other hand, it follows that
\begin{equation}
\label{eq:4.24}
\left|\int_0^\infty G_D(t,x,y)\varphi(y)\,\dee y\right|\le(4\pi t)^{-\frac{1}{2}}\|\varphi\|_{L^1}\to 0\quad\mbox{as}\quad t\to\infty.
\end{equation}
Furthermore, for any $\epsilon\in(0,1)$ and $t>0$, 
we have   
\begin{equation}
\label{eq:4.25}
\begin{split}
 & \left|\int_0^t \partial_yG_D(t-\tau,x,0)u(\tau,0)\,\dee\tau\right|
=\left|\int_0^t \partial_yG_D(\tau,x,0)u(t-\tau,0)\,\dee\tau\right|\\
 & \le\sup_{\tau\in(0,\infty)}u(\tau,0) \left|\int_{(1-\epsilon)t}^t\partial_yG_D(\tau,x,0)\,\dee\tau\right|
+\sup_{\tau\in[\epsilon t,\infty)}u(\tau,0)\left|\int_0^{(1-\epsilon)t} \partial_yG_D(\tau,x,0)\,\dee\tau\right|.
\end{split}
\end{equation}
Since 
$$
\partial_yG_D(\tau,x,0)=(4\pi \tau)^{-\frac{1}{2}}\frac{x}{\tau}\exp\left(-\frac{|x|^2}{4\tau}\right),
$$
we obtain
\begin{align*}
 & \int_{(1-\epsilon)t}^t\partial_yG_D(\tau,x,0)\,\dee\tau
\le C\int_{(1-\epsilon)t}^t \tau^{-1}\,\dee\tau=C\log\frac{1}{1-\epsilon},\\
 & \int_0^\infty \partial_yG_D(\tau,x,0)\,\dee\tau=\frac{1}{2\sqrt{\pi}}\int_0^\infty 2\xi^{-\frac{1}{2}}e^{-\xi}\,\dee\xi
=\frac{1}{\sqrt{\pi}}\Gamma\left(\frac{1}{2}\right)=1,
\end{align*}
for $x\in(0,\infty)$.
These estimates, together with Lemma~\ref{Lemma:4.1} and \eqref{eq:4.25}, imply that 
\begin{equation*}
\begin{split}
 & \limsup_{t\to\infty}\sup_{x\in(0,\infty)}\left|\int_0^t \partial_yG_D(t-\tau,x,0)u(\tau,0)\,\dee\tau\right|\\
 & \le C\log\frac{1}{1-\epsilon}\sup_{\tau\in(0,\infty)}u(\tau,0) +\limsup_{t\to\infty}\sup_{\tau\in[\epsilon t,\infty)}u(\tau,0)\le C\log\frac{1}{1-\epsilon}.
\end{split}
\end{equation*}
Consequently, letting $\epsilon\to 0$, we obtain 
\begin{equation}
\label{eq:4.26}
\lim_{t\to\infty}\sup_{x\in(0,\infty)}\left|\int_0^t \partial_yG_D(t-\tau,x,0)u(\tau,0)\,\dee\tau\right|=0.
\end{equation}
Combining \eqref{eq:4.24} and \eqref{eq:4.26}, we conclude that $\lim_{t\to\infty}\|z(t)\|_{L^\infty}=0$. 
Hence, we obtain $\lim_{t\to\infty}\|u(t)\|_{L^\infty}=0$, and Proposition~\ref{Proposition:4.1} follows. 
$\Box$\vspace{5pt}

By Proposition~\ref{Proposition:4.1}, for any global-in-time solution~$(s,u)$ to problem~\eqref{SP}, 
we have
$$
\sup_{t>0}\|u(t)\|_{L^\infty}<\infty.
$$
In the following proposition, we apply arguments similar to those in the proof of Lemma~\ref{Lemma:4.1} to 
obtain uniform $L^\infty$ estimates of global-in-time solutions to problem~\eqref{SP}.
\begin{proposition}
\label{Proposition:4.2}
For any $M>0$, there exists $C>0$ such that 
$$
\sup_{t\ge0}\|u(t)\|_{L^\infty}\le C
$$
for all global-in-time solution $(s,u)$ to problem~\eqref{SP} with initial data $(s_0,\varphi)$ satisfying 
$$
s_0+s_0^{-1}+\|\varphi\|_{{\rm Lip}([0,s_0])}\le M.
$$
\end{proposition}
{\bf Proof.}
The proof is by contradiction. 
Assume that there exist $M>1$ and 
a sequence of global-in-time solutions $\{(s_n,u_n)\}$ to problem~\eqref{SP} 
such that 
\begin{equation}
\label{eq:4.27}
\sup_n s_n(0)+\sup_n s_n(0)^{-1}\le M,\quad \sup_n\|u_n(0)\|_{{\rm Lip}([0,s_n(0)])}\le M,\quad \lim_{n\to\infty}\sup_{t\ge0}\|u_n(t)\|_{L^\infty}=\infty. 
\end{equation}
By Proposition~\ref{Proposition:2.1}, there exists $\delta>0$ such that 
\begin{equation}
\label{eq:4.28}
\sup_n \sup_{t\in[0,\delta]}\|u_n(t)\|_{L^\infty}\le 2M. 
\end{equation}
Then, by Proposition~\ref{Proposition:4.1} and the maximum principle, 
for any sufficiently large $n$, there exists a monotonically increasing sequence $\{t_n\}\subset(\delta,\infty)$ such that 
\begin{equation}
\label{eq:4.29}
\sigma_n:=\sup_{t\ge0}\|u_n(t)\|_{L^\infty}=\|u_n(t_n)\|_{L^\infty}=u_n(t_n,0)>2M. 
\end{equation}
For each $n=1,2,\dots$, similarly to \eqref{eq:4.4},
we define 
\begin{equation*}
\lambda_n:=\sigma_n^{-(p-1)},\quad \tau_n:=\lambda_n^{-2}t_n,\quad
s_n(\tau):=\lambda_n^{-1}s(t_n+\lambda_n^2\tau)\quad\mbox{for}\quad\tau\in[-\tau_n,0].
\end{equation*}
It follows from \eqref{eq:4.27} that $\lambda_n\to 0$ as $n\to\infty$.
Furthermore, it also from \eqref{eq:4.28} and \eqref{eq:4.29} that $t_n\ge\delta$ for $n=1,2,\dots$. 
Then
$$
\tau_n\to \infty\quad\mbox{as}\quad n\to\infty.
$$
Let $D_n$ and $D_n'$ be defined as in \eqref{eq:4.5}, and define 
$$
v_n(\tau,y)=\lambda_n^{\frac{1}{p-1}}u_n(t_n+\lambda_n^2\tau,\lambda_n y),\quad (\tau,y)\in D'_n.
$$
Then, it follows from \eqref{eq:4.29} that 
$$
v_n(0,0)=1,
\qquad
v_n\le 1\quad\mbox{in}\quad D'_n,\quad\mbox{for}\quad n=1,2,\dots. 
$$
Furthermore, similarly to \eqref{eq:4.7}, 
$v_n$ satisfies
\begin{equation*}
\left\{
\begin{array}{ll}
 \partial_{\tau}v_n=\partial_{y}^2v_n, & \quad (\tau,y)\in D_n,\vspace{7pt}\\
 -\partial_y v_n(\tau,0)=v_n(\tau,0)^p, & \quad \tau\in(-\tau_n,0],\vspace{7pt}\\
v_n(\tau,s_n(\tau))=0, & \quad \tau\in(-\tau_n,0],\vspace{3pt}\\
-\partial_y v_n(\tau,s_n(\tau))=\lambda_n^{\frac{1}{p-1}}s_n'(\tau), & \quad\tau\in(-\tau_n,0],
\end{array}
\right.
\end{equation*}
for $n=1,2\dots$. 
Then, applying arguments similar to those in the proof of Lemma~\ref{Lemma:4.1}, 
we arrive at a contradiction.
Thus, Proposition~\ref{Proposition:4.2} follows. 
$\Box$\vspace{5pt}

Furthermore, 
applying arguments similar to those in the proof of Lemma~\ref{Lemma:4.1} again, 
we prove the energy blow-up of solutions at the blow-up time.
\begin{proposition}
\label{Proposition:4.3}
Let $(s,u)$ be a solution to problem~\eqref{SP} with initial data $(s_0,\varphi) \in W$, 
and let $T_m$ denote its maximal existence time. 
If $T_m<\infty$, then 
$$
\lim_{t\nearrow T_m}E(s(t),u(t))=-\infty. 
$$
\end{proposition}
{\bf Proof.}
The proof is by contradiction. 
Assume that 
$$
\liminf_{t\nearrow T_m}E(s(t),u(t))>-\infty. 
$$
Then it follows from Lemma~\ref{Lemma:2.1} that 
\begin{equation}
\label{eq:4.30}
\int_{t_*}^{T_m}\int_0^{s(t)} |\partial_t u|^2\,\dee x\,\dee t+\frac{1}{2}\int_{t_*}^{T_m} (s')^3\,\dee t
\le E(s(t_*),u(t_*))-\liminf_{t\nearrow T_m}E(s(t),u(s))<\infty
\end{equation}
for $t_*\in(0,T_m)$. 
Since the solution blows up at $t=T_m$, 
by Proposition~\ref{Proposition:2.4}, there exists a monotonically increasing sequence $\{t_n\}\subset(0,T_m)$ such that 
\begin{equation}
\label{eq:4.31}
\lim_{n\to\infty}t_n=T_m,\quad
\lim_{n\to\infty}u(t_n,0)=\infty,\quad 
\sup_{t\in[t_0,t_n]} u(t,0)=u(t_n,0)\ge 1\quad\mbox{for $n=0,1,2,\dots$}.
\end{equation}
Let $M:=\|\varphi\|_{L^\infty}$.
By the maximum principle, 
for any sufficiently large $n$, we have
\begin{equation}
\label{eq:4.32}
\sigma_n:=\sup_{t\in(0,t_n]} \|u(t)\|_{L^\infty}=\|u(t_n)\|_{L^\infty}=u(t_n,0)>2M. 
\end{equation}
For each $n=1,2,\dots$, similarly to \eqref{eq:4.4},
we define 
\begin{equation*}
\lambda_n:=\sigma_n^{-(p-1)},\quad \tau_n:=\lambda_n^{-2}t_n,\quad
s_n(\tau):=\lambda_n^{-1}s(t_n+\lambda_n^2\tau)\quad\mbox{for}\quad\tau\in[-\tau_n,0].
\end{equation*}
It follows from \eqref{eq:4.27} that $\lambda_n\to 0$ as $n\to\infty$, which, together with \eqref{eq:4.31}, implies that 
$\tau_n\to \infty$ as $n\to\infty$. 
Let $D_n$ and $D_n'$ be defined as in \eqref{eq:4.5}, and define 
$$
v_n(\tau,y)=\lambda_n^{\frac{1}{p-1}}u(t_n+\lambda_n^2\tau,\lambda_n y),\quad (\tau,y)\in D'_n.
$$
Then it follows from \eqref{eq:4.32} that 
$$
v_n(0,0)=1,
\qquad
v_n\le 1\quad\mbox{in}\quad D'_n,\quad\mbox{for}\quad n=1,2,\dots. 
$$
Furthermore, similarly to \eqref{eq:4.7}, 
$v_n$ satisfies
\begin{equation*}
\left\{
\begin{array}{ll}
 \partial_{\tau}v_n=\partial_{y}^2v_n, & \quad (\tau,y)\in D_n,\vspace{7pt}\\
 -\partial_y v_n(\tau,0)=v_n(\tau,0)^p, & \quad \tau\in(-\tau_n,0],\vspace{7pt}\\
v_n(\tau,s_n(\tau))=0, & \quad \tau\in(-\tau_n,0],\vspace{3pt}\\
-\partial_y v_n(\tau,s_n(\tau))=\lambda_n^{\frac{1}{p-1}}s_n'(\tau), & \quad\tau\in(-\tau_n,0],
\end{array}
\right.
\end{equation*}
for $n=1,2\dots$. 
Then, by virtue of \eqref{eq:4.30}, 
applying arguments similar to those in the proof of Lemma~\ref{Lemma:4.1}, 
we arrive at a contradiction.
Thus, Proposition~\ref{Proposition:4.3} follows.
$\Box$
\vspace{5pt}

At the end of this section, we apply the comparison principle to derive a sufficient condition 
ensuring that the solution exists globally in time and decays exponentially as $t\to\infty$.
\begin{proposition}
\label{Proposition:4.4}
There exists $\delta>0$ such that, for any $(s_0,\varphi)\in W$,
if 
\begin{equation}
\label{eq:4.33}
\|\varphi\|_{L^\infty}\le \delta\min\left\{1,s_0^{-\frac{1}{p-1}}\right\},
\end{equation}
then problem~\eqref{SP} admits a global-in-time solution~$(s,u)$ with initial data $(s_0,\varphi)$ satisfying
\begin{equation}
\label{eq:4.34}
\begin{split}
 & s_\infty:=\lim_{t\to\infty}s(t)<\infty,\\
 & \|u(t)\|_{L^\infty}=O(e^{-\alpha t})\quad\mbox{as $t\to\infty$ for some $\alpha>0$}. 
\end{split}
\end{equation}
\end{proposition}
{\bf Proof.}
Let $s_0\in(0,\infty)$, and set 
\begin{equation}
\label{eq:4.35}
\epsilon:=\min\left\{\frac{1}{24},2^{-\frac{p+2}{p-1}}s_0^{-\frac{1}{p-1}}\right\}. 
\end{equation}
Choose $\alpha>0$ such that
\begin{equation}
\label{eq:4.36}
\frac{1}{32s_0^2}\le\alpha\le\frac{1}{16s_0^2}.
\end{equation}
Define
$$
\sigma(t):= 2s_0(2 - e^{-\alpha t})\text{ for $t\in [0, \infty)$},
\quad
V(y):=-y^2-y+2\in[0,2]\text{ for  $y\in[0,1]$}.
$$
Note that 
\begin{equation}
\label{eq:4.37}
2s_0\le \sigma(t)\le 4s_0,\quad \sigma'(t)>0, \quad t\in [0, \infty).
\end{equation}
Define
$$
v(t,x) := \epsilon e^{-\alpha t} V\left( \frac{x}{\sigma(t)} \right),\quad (t,x)\in D_\sigma. 
$$
It follows from $V(1)=0$ that
\begin{equation}
\label{eq:4.38}
v(t,\sigma(t))=0,\quad t\in[0,\infty).
\end{equation}
Since $V'(1)=-3$, we deduce from \eqref{eq:4.36} and \eqref{eq:4.37} that 
\begin{align*}
 & \sigma'(t)=2s_0\alpha e^{-\alpha t}\ge \frac{1}{16s_0}e^{-\alpha t}\ge\frac{3}{2s_0}\epsilon e^{-\alpha t},\\
 & -v_x(t,\sigma(t))=-\epsilon e^{-\alpha t}\sigma(t)^{-1}V'(1)\le\frac{3}{2s_0}\epsilon e^{-\alpha t},
\end{align*}
for $t\in(0,\infty)$. Hence, 
\begin{equation}
\label{eq:4.39}
\sigma'(t)\ge -v_x(t,\sigma(t)),\quad  t\in[0,\infty).
\end{equation}
Furthermore, since $V\le 2$, $V'\le 0$, and $V''=-2$ in $[0,1]$, 
we observe from \eqref{eq:4.36} and \eqref{eq:4.37} that 
\begin{equation}
\label{eq:4.40}
\begin{split}
\partial_tv -\partial_x^2 v 
 & = \epsilon e^{-\alpha t} \left( -\alpha V\left( \frac{x}{\sigma(t)} \right) - x \frac{\sigma'(t)}{\sigma(t)^2} V'\left( \frac{x}{\sigma(t)} \right)
- \frac{1}{\sigma(t)^2} V''\left( \frac{x}{\sigma(t)} \right)\right)\\
 & \ge \epsilon e^{-\alpha t} \left(-2\alpha+\frac{2}{(4s_0)^2}\right)\ge 0
\end{split}
\end{equation}
for $(t,x)\in D_\sigma$.
On the other hand, it follows from \eqref{eq:4.37}, $V(0)=2$, and $V'(0)=-1$ that
$$
-\partial_xv(t,0)=-\epsilon e^{-\alpha t}\sigma(t)^{-1}V'(0)=\epsilon e^{-\alpha t}\sigma(t)^{-1}
\ge \frac{1}{4s_0}\epsilon e^{-\alpha t},\quad 
v(t,0)^p=\left(2\epsilon e^{-\alpha t}\right)^p,
$$
for $t\in(0,\infty)$. These, together with \eqref{eq:4.35}, imply that  
\begin{equation}
\label{eq:4.41}
-\partial_xv(t,0)\ge v(t,0)^p,\quad t\in(0,\infty).
\end{equation}

Define
$$
\delta:=\frac{5}{4}\min\left\{\frac{1}{24},2^{-\frac{p+2}{p-1}}\right\}. 
$$
Assume \eqref{eq:4.33}. Since $V'\le 0$ on $[0,1]$, it follows that
$$
\min_{x\in[0,s_0]}v(0,x)=\epsilon V\left(\frac{1}{2}\right)=\frac{5}{4}\epsilon\ge \delta\min\left\{1,s_0^{-\frac{1}{p-1}}\right\}\ge\|\varphi\|_{L^\infty}.
$$
By \eqref{eq:4.38}--\eqref{eq:4.41}, and using the comparison principle (see the proof of Proposition~\ref{Proposition:2.5} in the case of \eqref{eq:2.39}), 
we deduce that problem~\eqref{SP} admits a global-in-time solution $(s,u)$ satisfying
$$
s(t)\le\sigma(t)\quad\mbox{for $t\in(0,\infty)$},\qquad u(t,x)\le v(t,x)\le 2\epsilon e^{-\alpha t}\quad\mbox{for $(t,x)\in D_s$}. 
$$
Consequently,
$$
\lim_{t\to\infty}s(t)\le\lim_{t\to\infty}\sigma(t)=4s_0,\qquad \|u(t)\|_{L^\infty}\le\|v(t)\|_{L^\infty}\le 2\epsilon e^{-\alpha t}\quad\mbox{for $t\in(0,\infty)$}.
$$
Hence, \eqref{eq:4.34} holds, and Proposition~\ref{Proposition:4.4} follows.
$\Box$
\vspace{5pt}

Combining Propositions~\ref{Proposition:2.6} and \ref{Proposition:4.4}, 
we obtain the stability of global-in-time solutions with exponential decay. 
\begin{proposition}
\label{Proposition:4.5}
Let $(s,u)$ be a global-in-time solution to problem~\eqref{SP} with initial data $(s_0,\varphi)\in W$.
Assume that $(s,u)$ satisfies \eqref{eq:4.34}.
Then there exists $\eta>0$ such that, 
for any $(\sigma_0,\psi)\in W$, 
if 
\begin{equation}
\label{eq:4.42}
|\sigma_0-s_0|+\|\psi-\varphi\|_{L^\infty}<\eta,
\end{equation}
then problem~\eqref{SP} admits a global-in-time solution $(\sigma,v)$ with initial data $(\sigma_0,\psi)$ such that 
\begin{equation}
\label{eq:4.43}
\begin{split}
 & \sigma_\infty:=\lim_{t\to\infty}\sigma(t)<\infty,\\
 & \|v(t)\|_{L^\infty}=O(e^{-\bar{\alpha} t})\quad\mbox{as $t\to\infty$ for some $\bar{\alpha}>0$}. 
\end{split}
\end{equation}
\end{proposition}
{\bf Proof.}
Let $(s,u)$ be a global-in-time solution to problem~\eqref{SP} with initial data $(s_0,\varphi)\in W$ 
and satisfy \eqref{eq:4.34}. Then we find $t_0>0$ such that 
$$
\|u(t_0)\|_{L^\infty}\le\frac{\delta}{2}\min\left\{1,s(t_0)^{-\frac{1}{p-1}}\right\},
$$
where $\delta$ is as in Proposition~\ref{Proposition:4.4}. 
Let $T:=t_0+1$ and $\eta>0$ be sufficiently small. 
Let $(\sigma,v)$ be a solution to problem~\eqref{SP} with initial data $(\sigma_0,\psi)\in W$, 
and assume \eqref{eq:4.42}. 
By Proposition~\ref{Proposition:2.6}, 
taking $\eta>0$ sufficiently small if necessary, 
we deduce that the solution~$(\sigma,v)$ exists on $(0,T]$ and 
\[
\|v(t_0)\|_{L^\infty}\le \delta\min\left\{1,\sigma(t_0)^{-\frac{1}{p-1}}\right\}.
\] 
Then, applying Propositions~\ref{Proposition:2.3} and \ref{Proposition:4.4}, 
we conclude that $(\sigma,v)$ is a global-in-time solution to problem~\eqref{SP} satisfying \eqref{eq:4.43}. 
Thus, Proposition~\ref{Proposition:4.5} follows. 
$\Box$
%%%%%%%%%%%%%%%%%%%%%%%%%%%%%%
%%%%%%%%%%%%%%%%%%%%%%%%%%%%%%
\section{Proofs of Theorems~\ref{Theorem:1.1} and \ref{Theorem:1.2}}
%%%%%%%%%%%%%%%%%%%%%%%%%%%%%%
%%%%%%%%%%%%%%%%%%%%%%%%%%%%%%
We complete the proofs of Theorems~\ref{Theorem:1.1} and \ref{Theorem:1.2}.
\vspace{5pt}
\newline
{\bf Proof of Theorem~\ref{Theorem:1.1}.}
Assertion~(1) follows from Proposition~\ref{Proposition:3.1}.
Assertion~(2) follows from Proposition~\ref{Proposition:2.4}, \ref{Proposition:3.3}, and \ref{Proposition:4.3}. 
Thus, Theorem~\ref{Theorem:1.1} follows. 
$\Box$
\vspace{5pt}
\newline
{\bf Proof of Theorem~\ref{Theorem:1.2}.}
Let $(s_0,\varphi)\in W$. 
For any $\lambda>0$, let $(s_\lambda,u_\lambda)$ be 
the solution to problem~\eqref{SP} with initial data $(s_0,\lambda\varphi)\in W$, 
and let $T_\lambda$ denote its maximal existence time. 

Let $\Lambda_*$ be the set of all $\lambda>0$ such that
$$
T_\lambda=\infty,\quad \lim_{t\to\infty}s_\lambda(t)<\infty,
\quad 
\|u_\lambda(t)\|_{L^\infty}=O(e^{-\alpha t})\mbox{ as $t\to\infty$ for some $\alpha>0$}.
$$
By Proposition~\ref{Proposition:4.4},  
there exists $\eta>0$ such that $(0,\eta)\subset\Lambda_*$. 
Furthermore, Proposition~\ref{Proposition:2.5} implies that 
if $\lambda\in\Lambda_*$, then $(0,\lambda]\subset\Lambda_*$. 
Then, by Proposition~\ref{Proposition:4.5}, we obtain
$\Lambda_*=(0,\lambda_*)$, where $\lambda_*:=\sup\Lambda_*$. 

We show that the set 
$$
\Lambda^*:=\{\lambda>0\,:\,T_\lambda=\infty\}
$$
is bounded. 
Consider the case where $\varphi(0)>0$. 
Let $\epsilon\in(0,1/4)$. 
Let $\zeta$ be a smooth nonnegative function on $[0,\infty)$ such that 
\begin{equation}
\label{eq:5.1}
\left\{
\begin{array}{ll}
\zeta(x)=1-x\quad & \mbox{for $x\in[0,\epsilon]$},\vspace{3pt}\\
|\partial_x\zeta(x)|\le 2\quad & \mbox{for $x\in[\epsilon,2\epsilon]$},\vspace{3pt}\\
\zeta(x)=1-2\epsilon & \mbox{for $x\in[2\epsilon,1]$},\vspace{3pt}\\
|\partial_x\zeta(x)|\le 2\epsilon & \mbox{for $x\in[1,1+\epsilon^{-1}]$},\vspace{3pt}\\
\zeta(x)=0 & \mbox{for $x\in[1+\epsilon^{-1},\infty$)}.
\end{array}
\right.
\end{equation}
Taking sufficiently small $\epsilon\in(0,1/4)$ if necessary, 
we obtain 
\begin{equation}
\label{eq:5.2}
\frac{1}{2}\int_0^\infty|\partial_x\zeta|^2\,\dee x-\frac{1}{p+1}\zeta(0)^{p+1}
\le \frac{1}{2}\left(5\epsilon+(2\epsilon)^2\epsilon^{-1}\right)-\frac{1}{p+1}=\frac{9}{2}\epsilon-\frac{1}{p+1}<0.
\end{equation}
Let $k>0$, and define
\begin{equation}
\label{eq:5.3}
\zeta_k(x):=k\zeta(k^{p-1}x),\quad x\in[0,\infty).
\end{equation}
Then
\begin{equation}
\label{eq:5.4}
-\partial_x\zeta_k(0)=k^p=\zeta_k(0)^p,
\end{equation}
and furthermore, it follows from \eqref{eq:5.2} that 
\begin{equation}
\label{eq:5.5}
\frac{1}{2}\int_0^\infty|\partial_x\zeta_k|^2\,\dee x-\frac{1}{p+1}\zeta_k(0)^{p+1}
=k^{p+1}\left(\frac{1}{2}\int_0^\infty|\partial_x\zeta|^2\,\dee x-\frac{1}{p+1}\zeta(0)^{p+1}\right)<0.
\end{equation}
On the other hand, 
since $\varphi$ is continuous on $[0,s_0]$ and $\varphi(0)>0$, 
there exists $\delta\in(0,s_0)$ such that $\varphi(x)\ge\varphi(0)/2$ for $x\in[0,\delta]$. 
We take sufficiently large $k>0$ so that 
$$
k^{p-1}\delta > 1+\epsilon^{-1}.
$$
Then, by \eqref{eq:5.1} and \eqref{eq:5.3}, 
taking sufficiently large $\lambda>0$, we obtain 
\begin{equation}
\label{eq:5.6}
\lambda\varphi(x)\ge\frac{\lambda}{2}\varphi(0)\ge k\ge\zeta_k(x)\quad\mbox{for $x\in[0,\delta]$},\qquad
\zeta_k=0\quad\mbox{for $x\in[\delta,s_0]$}.
\end{equation}
Furthermore, 
by Proposition~\ref{Proposition:3.1} and parabolic regularity theorems, together with \eqref{eq:5.4} and \eqref{eq:5.6},
problem~\eqref{SP} with initial data $(s_0,\zeta_k)\in W$ 
admits a solution $(\sigma_k,v_k)$ on $(0,T]$ for some $T\in(0,\infty)$ such that 
$$
\partial_x v_k\in C(D^{**}_{\sigma_k}(T)). 
$$  
Then, by \eqref{eq:5.5} and \eqref{eq:5.6}, 
we apply Propositions~\ref{Proposition:2.5} and \ref{Proposition:3.1} to obtain 
$T_\lambda\le T_{v_k}<\infty$, 
where $T_{v_k}$ is the maximal existence time of the solution~$(\sigma_k,v_k)$. 
This implies the boundedness of $\Lambda^*$ when $\varphi(0)>0$. 

Consider the case where $\varphi(0)=0$. 
Let $w$ be a solution to problem 
 $$
 \left\{
 \begin{array}{ll}
 \partial_t w=\partial_x^2 w,  & (t,x)\in(0,\infty)\times(0,s_0),\vspace{3pt}\\
 \partial_x w(t,0)=0, & t\in(0,\infty),\vspace{3pt}\\
 w(t,s_0)=0,& t\in(0,\infty),\vspace{3pt}\\
 w(0,x)=\varphi(x), & x\in[0,s_0].
 \end{array}
 \right.
 $$
 Let $\lambda\in\Lambda^*$. 
Since $\varphi\not\equiv 0$ and $\varphi\ge 0$ on $[0,s_0]$, 
by  the comparison principle, the maximum principle, and Hopf's lemma, we have 
 $$
 u_\lambda(t,x)\ge \lambda w(t,x)>0,\quad (t,x)\in(0,\infty)\times[0,s_0].
 $$
Since $\lambda\in\Lambda_*$, 
by Proposition~\ref{Proposition:2.5}, 
there exists a global-in-time solution $(\sigma_\lambda,W_\lambda)$ to problem~\eqref{SP} with initial data $(s_0,\lambda w(1))$ such that 
$$
u_\lambda(t+1,x)\ge W_\lambda(t,x),\quad (t,x)\in D_{\sigma_\lambda}.
$$
On the other hand, since $w(1,0)>0$, 
for any sufficiently large $\lambda$, $W_\lambda$ blows up in finite time 
which contradicts $\lambda\in\Lambda^*$ if $\lambda$ is sufficiently large. 
Thus, $\Lambda^*$ is bounded in the case where $\varphi(0)=0$. 
Hence, $\Lambda^*$ is bounded for $(s_0,\varphi)\in W$. 
Then, by Propositions~\ref{Proposition:2.6} and \ref{Proposition:4.2}, we obtain 
\begin{equation}
\label{eq:5.7}
0<\lambda_*\le \lambda^*:=\sup\Lambda^*\in\Lambda^*,
\qquad
\Lambda^*\setminus\Lambda_*=[\lambda_*,\lambda^*].
\end{equation}
Furthermore, by Proposition~\ref{Proposition:4.1}, we have 
\begin{equation}
\label{eq:5.8}
\lim_{t\to\infty}\|u_{\lambda}(t)\|_{L^\infty}=0,\quad \lambda\in\Lambda^*.
\end{equation}

Next, we show that $\lim_{t\to\infty}s_{\lambda_*}(t)=\infty$. 
It follows from \eqref{eq:5.7} and \eqref{eq:5.8} that
$$
\lim_{t\to\infty}\|u_{\lambda_*}(t)\|_{L^\infty}=0.
$$
Assume, for contradiction, that $\lim_{t\to\infty}s_{\lambda_*}(t)<\infty$. 
Then there exists $t_*>0$ such that 
$$
\|u_{\lambda_*}(t_*)\|_{L^\infty}\le \frac{1}{2}\delta\min\left\{1,s_{\lambda_*}(t_*)^{-\frac{1}{p-1}}\right\},
$$
where $\delta$ is given in Proposition~\ref{Proposition:4.4}. 
By Proposition~\ref{Proposition:2.6}, there exists $\epsilon>0$ such that
$$
\|u_{\lambda}(t_*)\|_{L^\infty}\le \delta\min\left\{1,s_{\lambda}(t_*)^{-\frac{1}{p-1}}\right\}
\quad\mbox{for}\quad \lambda\in[\lambda_*,\lambda_*+\epsilon).
$$
Then Proposition~\ref{Proposition:4.4} implies that $[\lambda_*,\lambda_*+\epsilon)\subset\Lambda_*$, 
which contradicts the definition of $\lambda_*$. 
Hence, $\lim_{t\to\infty}s_{\lambda_*}(t)=\infty$. 
Furthermore, Proposition~\ref{Proposition:2.5} implies that 
\begin{equation}
\label{eq:5.9}
\lim_{t\to\infty}s_\lambda(t)\ge\lim_{t\to\infty}s_{\lambda_*}(t)=\infty,\quad \lambda\in\Lambda^*\setminus\Lambda_*. 
\end{equation}

Let $\lambda\in[\lambda_*,\lambda^*]$.
It remains to prove that
\begin{align}
\label{eq:5.10}
 & \liminf_{t\to\infty}s_\lambda(t)^{\frac{1}{p-1}}\|u_\lambda(t)\|_{L^\infty(0,s_\lambda(t))}>0,\\
\label{eq:5.11}
 & s_\lambda(t)=(1+o(1))\int_0^t u_\lambda(t,0)^p\,\dee\tau=o\left(t^{\frac{1}{2}}\right)\mbox{ as $t\to\infty$}.
\end{align}
By Proposition~\ref{Proposition:4.4} and \eqref{eq:5.9}, we have 
\begin{equation}
\label{eq:5.12}
\|u_\lambda(t)\|_{L^\infty}>\delta\min\left\{1,s_\lambda(t)^{-\frac{1}{p-1}}\right\}
=\delta s_\lambda(t)^{-\frac{1}{p-1}}
\end{equation}
for sufficiently large $t>0$. This implies \eqref{eq:5.10}. 
Furthermore, by \eqref{eq:2.9}, \eqref{eq:5.8}, and \eqref{eq:5.9}, we have 
\begin{align*}
s_\lambda(t) & =s_0+\lambda\|\varphi\|_{L^1}-\|u_\lambda(t)\|_{L^1}+\int_0^t u_\lambda(\tau,0)^p\,\dee\tau\\
 & =s_0+\lambda\|\varphi\|_{L^1}+o(s_\lambda(t))+\int_0^t u_\lambda(\tau,0)^p\,\dee\tau
\end{align*}
as $t\to\infty$. Hence,   
\begin{equation}
\label{eq:5.13}
s_\lambda(t)=(1+o(1))\int_0^t u_\lambda(\tau,0)^p\,\dee\tau
\end{equation}
as $t\to\infty$. 

Let $T>0$. By \eqref{eq:5.8}, there exists $t_T\in[T,\infty)$ such that 
\begin{equation}
\label{eq:5.14}
\sup_{\tau\in[T,\infty)}u_\lambda(\tau,0)=u_\lambda(t_T,0)=\sup_{\tau\in[t_T,\infty)}u_\lambda(\tau,0).
\end{equation}
Let $(\sigma_T,z_T)$ be a solution to the Stefan problem for the heat equation with an inhomogeneous Dirichlet boundary condition, 
\begin{equation*}
\left\{
\begin{array}{ll}
\partial_t z=\partial_x^2 z, & (t,x)\in D_\sigma,\vspace{3pt}\\
z(t,0)=u_\lambda(t_T,0),\qquad & t\in(0,\infty), \vspace{3pt}\\
z(t,\sigma(t))=0, &  t\in (0,\infty),\vspace{3pt}\\
\sigma'(t)=-\partial_xz(t,\sigma(t)), &  t\in (0,\infty), \vspace{3pt}\\
z(0,x)=u_\lambda(t_T,x), & x\in[0,s_\lambda(t_T)],\vspace{3pt}\\
\sigma(0)=s_\lambda(t_T).
\end{array}
\right.
\end{equation*}
By \eqref{eq:5.14}, the comparison principle yields
$$
s_\lambda(t+t_T)\le \sigma_T(t)\mbox{ on $[0,\infty)$},\qquad u_\lambda(t+t_T,x)\le z_T(t,x)\mbox{ in $D_{s_{\lambda,T}}$},
$$
where $s_{\lambda,T}(t):=s_\lambda(t+t_T)$ for $t\ge 0$. 
On the other hand, by \cites{BHMS, M}, we see that 
$$
\lim_{t\to\infty}t^{-\frac{1}{2}}\sigma_T(t)=A_T,
$$
where $A_T$ is the unique solution $A$ to the equation
$$
u_\lambda(t_T,0)=\frac{A}{2}e^{\frac{A^2}{4}}\int_0^A e^{-\frac{\tau^2}{4}}\,\dee\tau. 
$$
Note that $A_T\to 0$ as $T\to\infty$, since $u_\lambda(t_T,0)\to 0$ as $T\to\infty$.
Then we obtain
\begin{equation}
\label{eq:5.15}
\limsup_{t\to\infty}t^{-\frac{1}{2}}s_\lambda(t)=\limsup_{T\to\infty}\limsup_{t\to\infty}t^{-\frac{1}{2}}s_{\lambda,T}(t)
\le\lim_{T\to\infty}\lim_{t\to\infty}t^{-\frac{1}{2}}\sigma_T(t)=\lim_{T\to\infty}A_T=0.
\end{equation}
This, together with \eqref{eq:5.13}, implies \eqref{eq:5.11}.
Hence, the proof of Theorem~\ref{Theorem:1.2} is complete.
$\Box$\vspace{5pt}
\newline
{\bf Proof of Theorem~\ref{Theorem:1.3}.} 
Set $w:=\partial_x u_\lambda$. It follows that 
\begin{equation}
\left\{
\begin{array}{ll}
\partial_t w=\partial_x^2 w,\quad & (t,x)\in D_{s_\lambda},\vspace{3pt}\\
w(t,0)\le 0, & t\in[0,\infty),\vspace{3pt}\\
w(t,s_\lambda(t))\le 0, & t\in[0,\infty),\\
w(0,x)\le 0, & x\in(0,s_0).
\end{array}
\right.
\end{equation}
Then the maximum principle implies that $w\le 0$ in $D_{s_{\lambda}}$, and hence 
$\partial_x u_{\lambda}\le 0$ in $D_{s_\lambda}$. 
Consequently, we obtain
\begin{equation}
\label{eq:5.17}
\|u_\lambda(t)\|_{L^\infty}=u_\lambda(t,0),\quad t>0.
\end{equation}

Let $\lambda\in[\lambda_*,\lambda^*]$. 
By \eqref{eq:5.10}, \eqref{eq:5.11}, and \eqref{eq:5.17}, 
there exists $T>0$ such that
\begin{equation}
\label{eq:5.18}
u_\lambda(t,0)=\|u_\lambda(t)\|_{L^\infty}\ge Cs_\lambda(t)^{-\frac{1}{p-1}}
\ge C\left(\int_0^t u_\lambda(\tau,0)^p\,\dee\tau\right)^{-\frac{1}{p-1}}
\end{equation}
for $t\in[T,\infty)$. 
Setting
$$
X(t):=\int_0^t u_\lambda(\tau,0)^p\,\dee\tau,\quad t\in[T,\infty), 
$$
by \eqref{eq:5.18}, we have 
$$
\left(X'(t)\right)^{\frac{1}{p}}\ge CX(t)^{-\frac{1}{p-1}},\quad t\in[T,\infty), 
$$
that is, 
$$
X(t)^{\frac{p}{p-1}}X'(t)\ge C,\quad t\in[T,\infty).
$$
This implies that 
$$
X(t)^{\frac{2p-1}{p-1}}-X(T)^{\frac{2p-1}{p-1}}\ge C(t-T),\quad t\in[T,\infty). 
$$
Consequently, we obtain 
$$
\int_0^t u_\lambda(\tau,0)^p\,\dee\tau=X(t)\ge Ct^{\frac{p-1}{2p-1}}
$$
for sufficiently large $t$, which, together with \eqref{eq:5.11}, yields 
$s_\lambda(t)\ge  Ct^{\frac{p-1}{2p-1}}$
for sufficiently large $t$. Thus, \eqref{eq:1.13}, and Theorem~\ref{Theorem:1.3} follows.
$\Box$
\vspace{8pt}

\noindent
{\bf Acknowledgment.}
K. I. was supported in part by JSPS KAKENHI Grant Number 25H00591. 
%%%%%%%%%%%%%%%%%%%%%%%%%%%%%%%%%%%%%%
%%%%%%%%%%%%    references    %%%%%%%%%%%%%%%%%%
%%%%%%%%%%%%%%%%%%%%%%%%%%%%%%%%%%%%%%

\end{document}